\input epsf
%\magnification=\magstep1
\documentstyle{amsppt}
\pagewidth{6.0truein}%\hcorrection{0.55in}
\pageheight{9truein}%\vcorrection{0.75in}
\TagsOnRight
\NoRunningHeads
\catcode`\@=11
\def\logo@{}
\footline={\ifnum\pageno>1 \hfil\folio\hfil\else\hfil\fi}
\topmatter
\title 
%Four factorization formulas of the same kind for plane partitions
%A unified set of four identities involving all ten symmetry classes of plane partitions
%A unified set of four identities for plane partitions
Four factorization formulas for plane partitions
\endtitle
\vskip-0.2in
\author Mihai Ciucu\endauthor
\thanks Research supported in part by NSF grant DMS-1101670.
\endthanks
\affil
  Department of Mathematics, Indiana University\\
  Bloomington, Indiana 47405
\endaffil
\abstract
All ten symmetry classes of plane partitions that fit in a given box are known to be enumerated by simple product formulas, but there is still no unified proof for all of them. Progress towards this goal can be made by establishing identities connecting the various symmetry classes. We present in this paper four such identities, involving all ten symmetry classes. We discuss their proofs and generalizations. The main result of this paper is to give a generalization of one of them, in the style of the identity presented in  ``A factorization theorem for rhombus tilings,'' M. Ciucu and C. Krattenthaler, arXiv:1403.3323.
\endabstract
\endtopmatter

\document

\def\mysec#1{\bigskip\centerline{\bf #1}\message{ * }\nopagebreak\bigskip\par}

\def\myref#1{\item"{[{\bf #1}]}"}

\def\epf{\hfill{$\square$}\smallpagebreak}

\def\cite#1{\relaxnext@
  \def\nextiii@##1,##2\end@{[{\bf##1},\,##2]}%
  \in@,{#1}\ifin@\def\next{\nextiii@#1\end@}\else
  \def\next{[{\bf#1}]}\fi\next}
\def\proclaimheadfont@{\smc}

\define\Z{{\Bbb Z}}

\define\M{\operatorname{M}}

%\define\M{\operatorname{M}}

\define\twoline#1#2{\line{\hfill{\smc #1}\hfill{\smc #2}\hfill}}
\define\twolinetwo#1#2{\line{{\smc #1}\hfill{\smc #2}}}
\define\twolinethree#1#2{\line{\phantom{poco}{\smc #1}\hfill{\smc #2}\phantom{poco}}}
\define\threeline#1#2#3{\line{\hfill{\smc #1}\hfill{\smc #2}\hfill{\smc #3}\hfill}}

\def\mypic#1{\epsffile{#1}}

%\def\epsfsize#1#2{0.36#1}
%\topinsert
%\centerline{\mypic{2-1.eps}}
%\centerline{{\smc Figure~2.1{\rm (a). $R_{(2,4,5),(2,4)}(2)$.}}}
%\endinsert

%\topinsert
%\twoline{\mypic{2-1b.eps}}{\mypic{2-1c.eps}}
%\twoline{Figure~2.1{\rm (b). $R_{\emptyset,(2,4)}(4)$.}}{Figure~2.1{\rm (c). 
%$R_{(2,4,5),\emptyset}(2)$.}}
%\endinsert

% ref nos
\define\AndTSSC{1}
\define\FT{2}
\define\ppone{3}
\define\csts{4}
\define\ranglep{5}
\define\fakt{6}
\define\DT{7}
\define\KupCSSC{8}
\define\MacM{9}
\define\StanPP{10}
\define\StemTS{11}

% eq nos
\define\eaa{1.1}
\define\eab{1.2}
\define\eac{1.3}
\define\ead{1.4}
\define\eae{1.5}
\define\eaf{1.6}
\define\eag{1.7}
\define\eah{1.8}
\define\eai{1.9}
\define\eaj{1.10}
\define\eak{1.11}
\define\eal{1.12}

\define\eba{2.1}
\define\ebb{2.2}
\define\ebc{2.3}
\define\ebd{2.4}

\define\eca{3.1}
\define\ecb{3.2}
\define\ecc{3.3}
\define\ecd{3.4}
\define\ece{3.5}
\define\ecf{3.6}
\define\ecg{3.7}
\define\ech{3.8}
\define\eci{3.9}
\define\ecj{3.10}
\define\eck{3.11}
\define\ecl{3.12}
\define\ecm{3.13}

%\define\edb{4.2}

% th nos

\define\tba{2.1}

% fig nos
\define\faa{1.1}
\define\fab{1.2}
\define\fac{1.3}

\define\fba{2.1}
\define\fbb{2.2}

\define\fca{3.1}
\define\fcb{3.2}
\define\fcc{3.3}
\define\fcd{3.4}
\define\fce{3.5}
\define\fcf{3.6}
\define\fcg{3.7}
\define\fch{3.8}
\define\fci{3.9}
\define\fcj{3.10}
\define\fck{3.11}
\define\fcl{3.12}

\vskip-0.05in
\mysec{1. Introduction}

All ten symmetry classes of plane partitions that fit in a given box are known to be enumerated by beautifully elegant simple product formulas (\cite{\StanPP} covers seven symmetry classes, while \cite{\AndTSSC}, \cite{\KupCSSC} and \cite{\StemTS} prove the remaining three). 

Using these formulas, one readily checks that for any non-negative integers $a$ and $b$ one has
$$
\spreadlines{2\jot}
\align
P(a,a,2b)&=S(a,a,2b)\,TC(a,a,2b)\tag\eaa\\
SC(a,a,2b)&=SSC(a,a,2b)^2\tag\eab\\
CS(2a,2a,2a)&=TS(2a,2a,2a)\,CSTC(2a,2a,2a)\tag\eac\\
CSSC(2a,2a,2a)&=TSSC(2a,2a,2a)^2,\tag\ead
\endalign
$$
where the notation is explained in Table 1.1 (which is adapted from \cite\StanPP); the arguments indicate the sides of the box containing the plane partitions\footnote{ We note that in each of these four equalities, the arguments are in the most general form for which all the involved quantities are defined.}. An interesting feature of this set of identities is that it features all ten symmetry classes of plane partitions.

\topinsert
$$
\hbox{\vbox{\offinterlineskip
  \def\spacing{\omit\vrule height4pt&&&&&&\cr}
  \halign{\vrule\strut#&&\ \hfil$#$\hfil\ &\vrule#\cr
  \noalign{\hrule}\spacing
  & &&\text{Notation}&&\text{Symmetry class}&\cr
  \spacing\noalign{\hrule}\spacing
  &1&&P&&\text{Any}&\cr\spacing
  &2&&S&&\text{Symmetric}&\cr\spacing
  &3&&CS&&\text{Cyclically symmetric}&\cr\spacing
  &4&&TS&&\text{Totally symmetric}&\cr\spacing
  &5&&SC&&\text{Self-complementary}&\cr\spacing
  &6&&TC&&\text{Transpose-complementary (complement=transpose)}&\cr\spacing
  &7&&SSC&&\text{Symmetric and self-complementary}&\cr\spacing
  &8&&CSTC&&\text{Cyclically symmetric and transpose-complementary}&\cr\spacing
  &9&&CSSC&&\text{Cyclically symmetric and self-complementary}&\cr\spacing
  &10&&TSSC&&\text{Totally symmetric and self-complementary}&\cr\spacing
  \noalign{\hrule}}}}
$$
\centerline{{\smc Table~1.1.} Notation for the ten symmetry classes of plane partitions.}
\endinsert

A unified way of viewing these four equalities is afforded by the well-known equivalence between plane partitions that fit in an $a\times b\times c$ box and lozenge tilings of a hexagon of side-lengths $a$, $b$, $c$, $a$, $b$, $c$ (in cyclic order) on the triangular lattice (see e.g. \cite{\DT}; for definiteness, we consider that the triangular lattice is drawn in the plane so that one family of parallel lattice lines is vertical). 

Furthermore, the ten symmetry classes of plane partitions correspond to the symmetry classes of lozenge tilings of the corresponding hexagon. More precisely, the generators of the group of symmetries of plane partitions, namely $(i)$ swapping the $x$ and $y$ coordinate axes, $(ii)$ cyclically shifting the axes, and $(iii)$ taking the complement of the diagram of the plane partition in the box enclosing it, correspond to the associated tiling being invariant under reflection across the vertical, rotation by $120$ degrees, and rotation by $180$ degrees, respectively. 

Denote by $H_{a,b,c}$ the hexagon of sides $a$, $b$, $c$, $a$, $b$, $c$ (in clockwise order, starting from the northwestern side) on the triangular lattice; write $H_a$ for $H_{a,a,a}$. Translated in the language of tilings, equations (\eaa)--(\ead) read
$$
\spreadlines{2\jot}
\align
\M(H_{a,a,2b})&=\M_{-}(H_{a,a,2b})\M_{|}(H_{a,a,2b})\tag\eae\\
\M_{r^3}(H_{a,a,2b})&=\left(\M_{r^3,|}(H_{a,a,2b})\right)^2\tag\eaf\\
\M_{r^2}(H_{2a})&=\M_{r^2,-}(H_{2a})\M_{r^2,|}(H_{2a})\tag\eag\\
\M_{r}(H_{2a})&=\left(\M_{r,|}(H_{2a})\right)^2,\tag\eah
\endalign
$$
where $\M(R)$ denotes the number of lozenge tilings of the lattice region $R$, $\M_{s_1,\dotsc,s_k}(R)$ denotes the number of these which are invariant under the symmetries $s_1,\dotsc,s_k$, \linebreak while $-$, $|$, and $r$ denote reflection across the horizontal, reflection across the vertical, and rotation by $60^\circ$, respectively.

This set of identities already looks quite uniform, but it can be made even more so by noticing that equations (\eaf)--(\eah) can be formally obtained from (\eae) by ``moding out'' by the action of the cyclic groups $\Z_2$, $\Z_3$ and $\Z_6$, respectively. 

To explain this, it will be useful to recall that the lozenge tilings of any region $R$ on the triangular lattice can naturally be identified with the perfect matchings of the planar dual graph of $R$, i.e. the graph whose vertices are the unit triangles in $R$, and whose edges connect vertices corresponding to unit triangles that share an edge (this will be a subgraph of the hexagonal lattice). In order to keep the notation simpler --- and in view of the identification mentioned above --- we denote the dual graph of a region by the same symbol as the region itself. Thus $R$ may denote a lattice region, or its planar dual graph, depending on the context. In the latter case, $\M(R)$ will denote the number of perfect matchings (often simply referred to as matchings) of $R$.

\topinsert
\centerline{\mypic{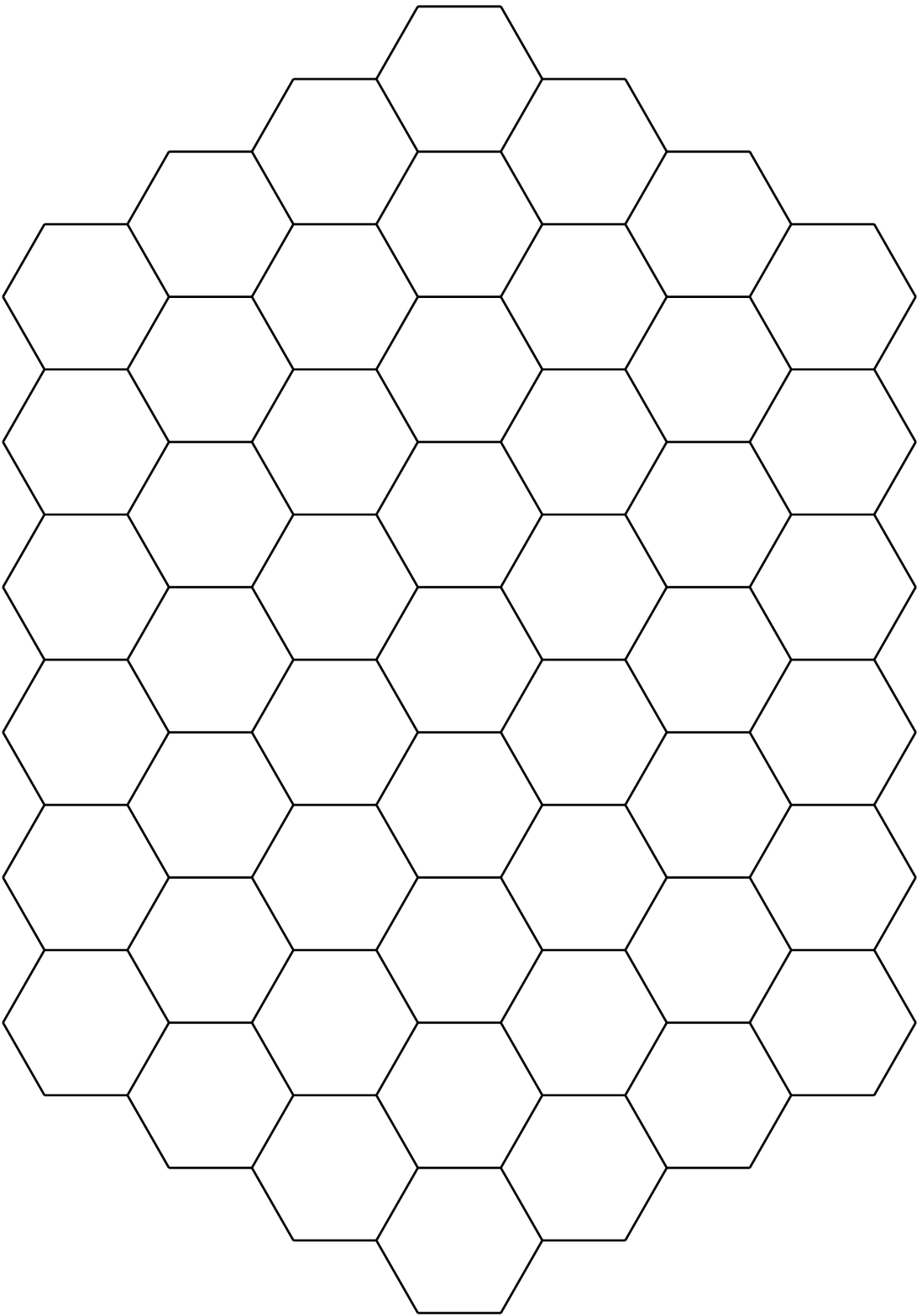}}
\twoline{\mypic{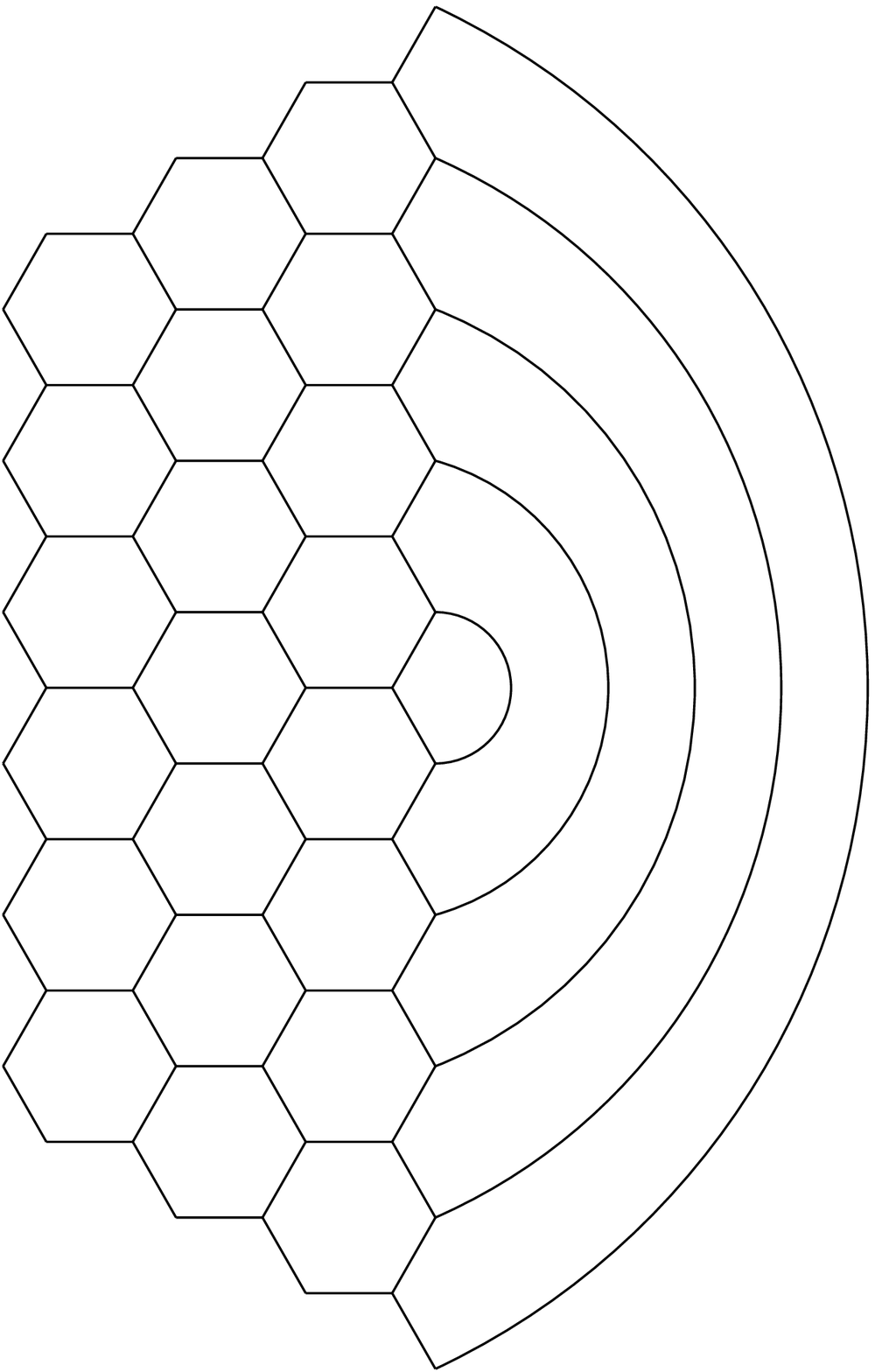}}{\mypic{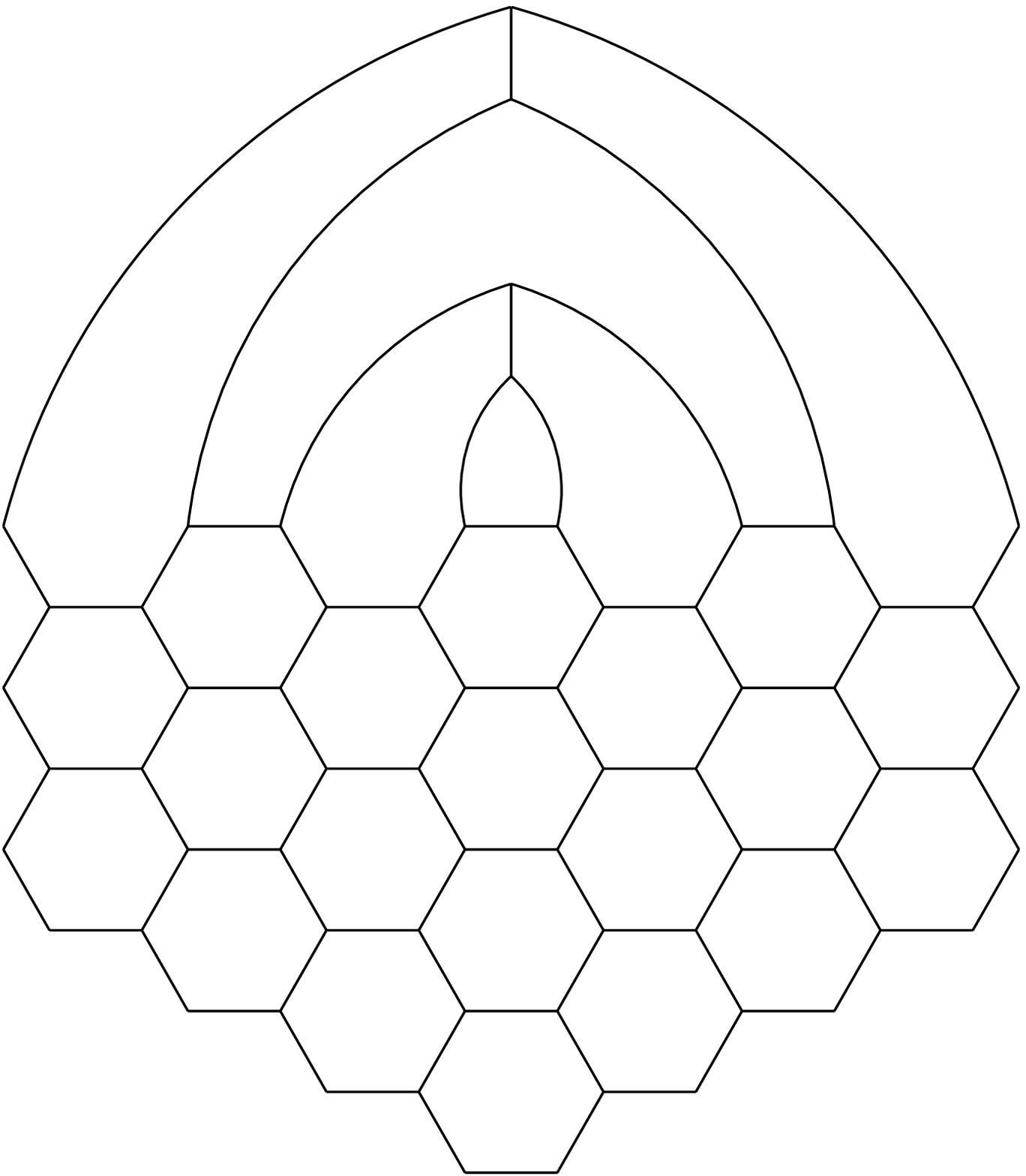}}
\medskip
\centerline{{\smc Figure~{\faa}. {\rm  The planar dual of $H_{4,4,6}$ and two drawings of its quotient by $\Z_2$.}}}
\endinsert

Consider the action of the group generated by the $180^\circ$ rotation around its center on the dual graph of $H_{a,a,2b}$, and denote by $H_{a,a,2b}/\Z_2$ the orbit graph (we will also refer to it as the quotient graph). Figure {\faa} shows the graph $H_{4,4,6}$ and illustrates two different ways of drawing its quotient in the plane: the first can be used to see what it means for a perfect matching of the quotient to be horizontally symmetric, and the second to see what it means to be vertically symmetric. Then, formally, equation (\eae) is mapped to
$$
\M(H_{a,a,2b}/\Z_2)=\M_{-}(H_{a,a,2b}/\Z_2)\M_{|}(H_{a,a,2b}/\Z_2);
$$
the remarkable thing is that
%(see Figure {\faa} for an illustration). Note that under this quotient, the vertical symmetry axis remains a symmetry axis of the resulting graph, while the horizontal symmetry axis is mapped to this vertical symmetry axis. Therefore, formally, equation (\eae) is mapped to
%, for even\footnote{ If $b$ is odd, it turns out that $\M_{|}(H_{a,a,b}/\Z_2)=0$. Indeed, by the paragraph following this footnote, vertically symmetric perfect matchings of $H_{a,a,b}/\Z_2$ correspond to perfect matching of $H_{a,a,b}$ that are invariant under both rotation by $180^\circ$ and reflection across the vertical, and are therefore necessarily also symmetric across the horizonta. However, a clear necesary condition for such perfect matchings to exist is that $b$ is even.} $b$, 
this is equivalent to (\eaf) --- and therefore a true equality!

To see this equivalence, note first that the isomorphism between the bottom two graphs in Figure {\faa} maps their symmetry axes to one another (this is evident if the graphs are embedded in the surface of a cone in such a way that all hexagonal faces are congruent). Due to this, the above equality is equivalent to
$$
\M(H_{a,a,2b}/\Z_2)=\M_{|}(H_{a,a,2b}/\Z_2)^2.
$$
To complete the argument, note that the perfect matchings of $H_{a,a,2b}/\Z_2$ can be identified with those perfect matchings of $H_{a,a,2b}$ which are invariant under rotation by $180^\circ$, and that under this bijection the vertically symmetric perfect matchings of $H_{a,a,2b}/\Z_2$ correspond to perfect matchings of $H_{a,a,2b}$ that are invariant under both rotation by $180^\circ$ and reflection across the vertical. The latter correspond precisely to symmetric and self-complementary plane partitions, according to the paragraph before equations (\eae)--(\eah).

\topinsert
\centerline{\mypic{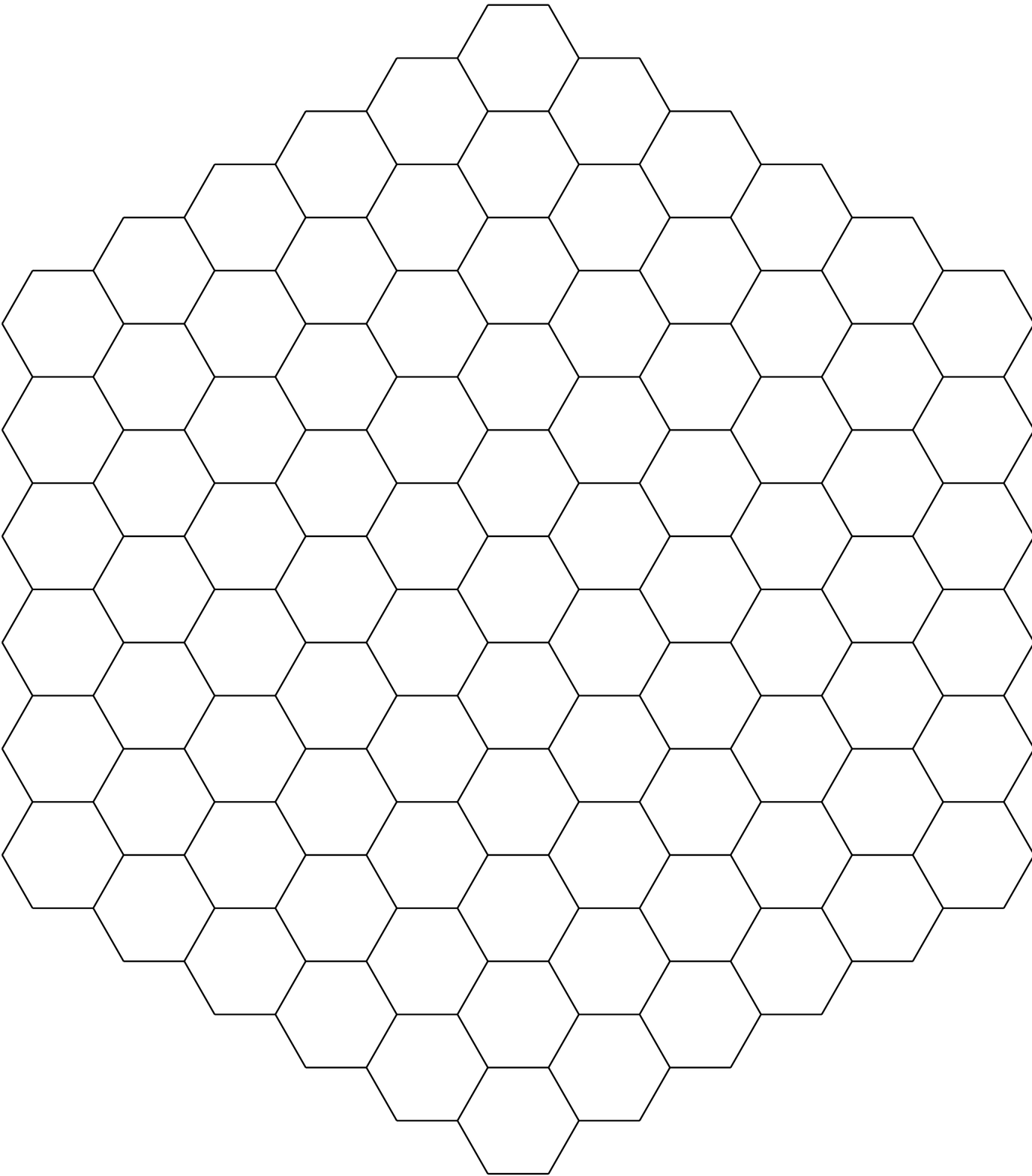}}
\twoline{\mypic{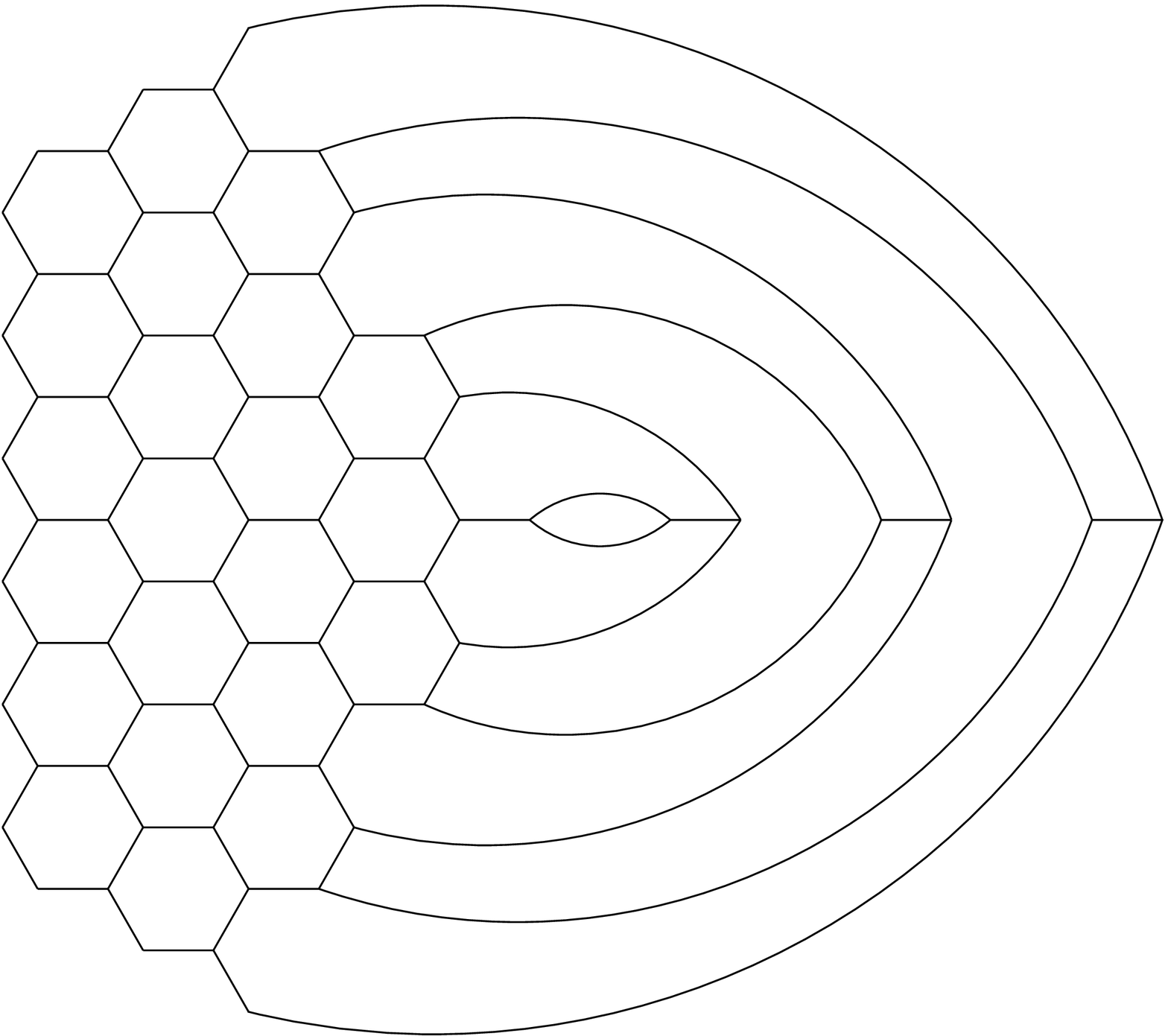}}{\mypic{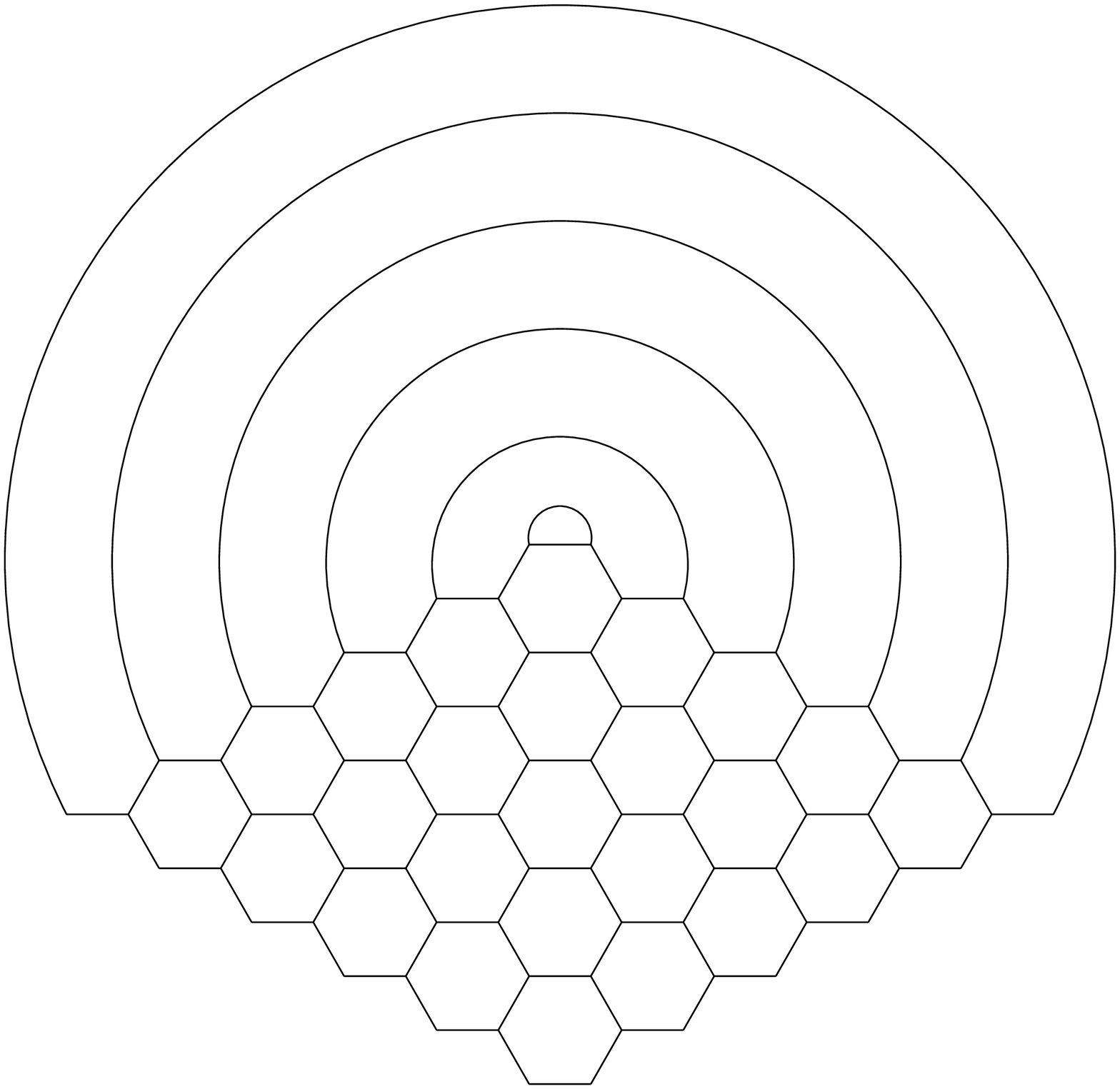}}
\medskip
\centerline{{\smc Figure~{\fab}. {\rm  The planar dual of $H_{6}$ and two drawings of its quotient by $\Z_3$.}}}
\endinsert

In the same way, one sees that (\eag) and (\eah) are obtained from (\eae) by moding out by $\Z_3$ and $\Z_6$, respectively.

More precisely, let $H_{2a}/\Z_3$ be the orbit (or quotient) graph of $H_{2a}$ under the action of the group generated by rotation by $120^\circ$ around its center. Figure {\fab} illustrates $H_6$ and two different drawings of its quotient by this rotation. The first drawing shows what it means for a matching of the quotient graph to be horizontally symmetric, and the second what it means to be vertically symmetric. Then, formally, when quotioning out by the $120^\circ$ rotation around the center, equation (\eae) is mapped to
$$
\M(H_{2a}/\Z_3)=\M_{-}(H_{2a}/\Z_3)\M_{|}(H_{2a}/\Z_3),
$$
%Call the three symmetry axes of $H_{2a}$ that are images of the horizontal symmetry axes under $120^\circ$ rotations red (they are indicated by solid lines in Figure {\fab}), and the remaining three symmetry axes (images of the vertical symmetry axis under rotations by $120^\circ$) blue (they are indicated by dashed lines in Figure~{\fab}). 
%Under the quotient, all three red symmetry axes of $H_{2a}$ are mapped to the horizontal, and all three blue axes to the vertical. Therefore, under this quotient map, (\eae) formally becomes
and again
%, for even\footnote{ A necessary condition for $\M_{-}(H_{a}/\Z_3)$ to be non-zero; see also footnote 1.} $a$, 
this turns out to be a true equality, this time due to the fact that it is equivalent to (\eag)! This equivalence holds because the perfect matchings of $H_{2a}/\Z_3$ are naturally identified with the perfect matchings of $H_{2a}$ that are invariant under rotation by $120^\circ$, and the subsets of the latter that are in addition symmetric across the horizontal (resp., vertical) correspond to perfect matchings of $H_{2a}/\Z_3$ that posses the additional symmetry.

\topinsert
\twoline{\mypic{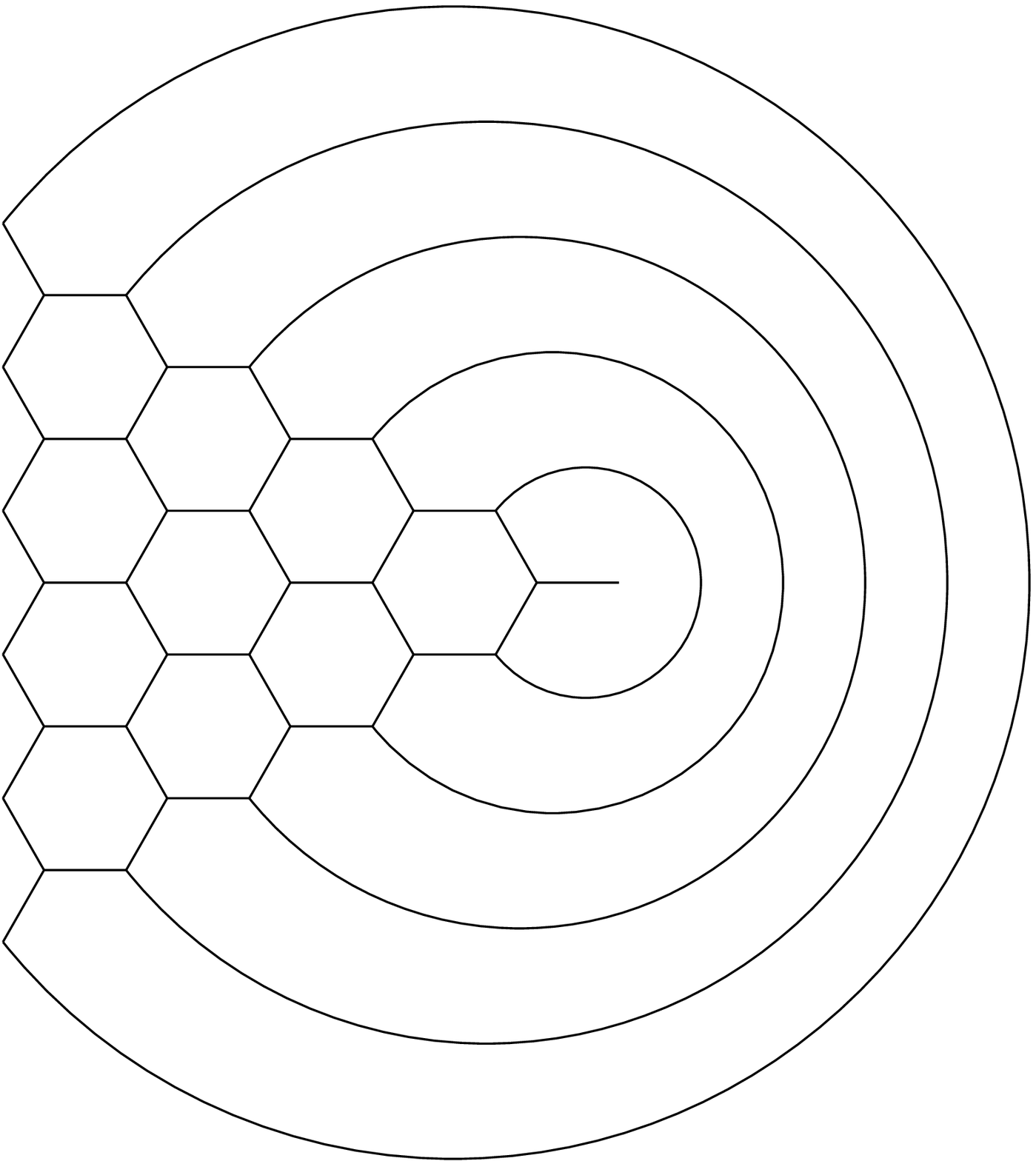}}{\mypic{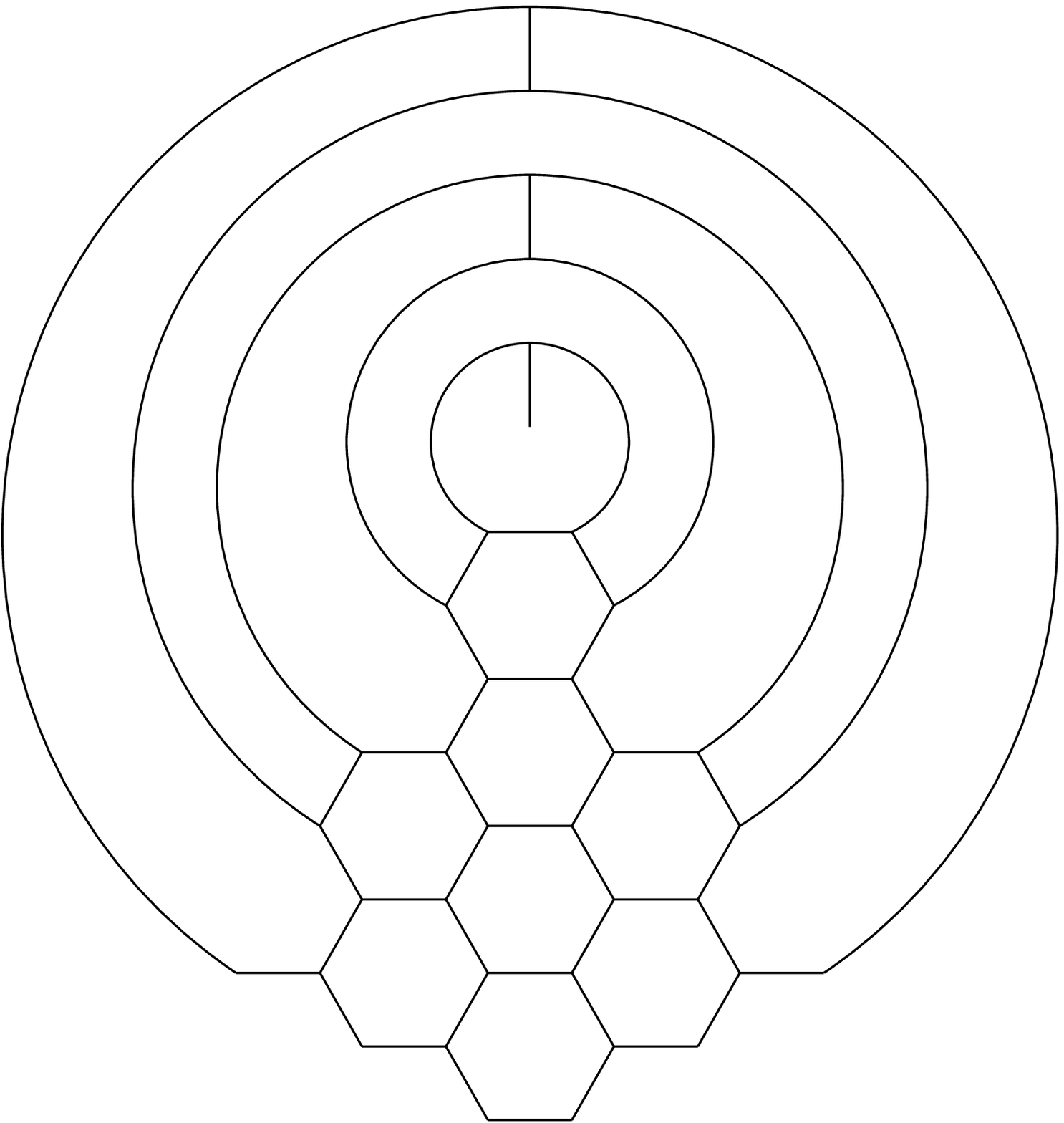}}
\medskip
\centerline{{\smc Figure~{\fac}. {\rm  Two drawings of the quotient of $H_6$ by $\Z_6$.}}}
\endinsert

In order to interpret (\ead) in the same style, consider the orbit graph $H_{2a}/\Z_6$ of the action of the group generated by the $60^\circ$ degree rotation around the center of $H_{2a}$. Figure~{\fac} shows two drawings of this quotient. The first can be used to define horizontally symmetric matchings of the quotient graphs, and the second to define vertically symmetric matchings. Then, formally, when quotioning out by the $60^\circ$ rotation around the center, equation (\eae) is mapped to
$$
\M(H_{2a}/\Z_6)=\M_{-}(H_{2a}/\Z_6)\M_{|}(H_{2a}/\Z_6),
$$
%
%
%Under this quotient, both the three red and the three blue symmetry axes of $H_{2a}$ are mapped to the vertical. Thus, under this quotient equation (\eae) formally becomes
%$$
%\M(H_{2a}/\Z_6)=\M_{|}(H_{2a}/\Z_6)^2,
%$$
which once more
%, for even\footnote{ In analogy to the previous two footnotes, $\M_{|}(H_{a}/\Z_6)=0$ for odd $a$.} values of $a$, 
turns out to be a true equality, due now to the fact that it is equivalent to~(\eah)! 

To see this, note that the two graphs in Figure {\fac} are isomorphic (as they represent the same quotient graph), and that their symmetry axes are mapped to one another via this isomorphism (this is apparent if the graph is drawn on a conical surface so that all hexagonal faces are congruent). It follows that the above equation is equivalent to 
$$
\M(H_{2a}/\Z_6)=\M_{|}(H_{2a}/\Z_6)^2.
$$
This in turn is equivalent to (\ead), because the perfect matchings of $H_{2a}/\Z_6$ can be identified with the perfect matchings of $H_{2a}$ that are invariant under rotation by  $60^\circ$, with the subset of the former that are also symmetric across the vertical corresponding to $60^\circ$ rotation invariant perfect matchings of $H_{2a}$ that are also invariant under reflection across the vertical.

Therefore, equations (\eae)--(\eah) can be rewritten as
$$
\spreadlines{2\jot}
\align
\M(H_{a,a,2b})&=\M_{-}(H_{a,a,2b})\M_{|}(H_{a,a,2b})\tag\eai\\
\M(H_{a,a,2b}/\Z_2)&=\M_{|}(H_{a,a,2b}/\Z_2)^2\tag\eaj\\
\M(H_{2a}/\Z_3)&=\M_{-}(H_{2a}/\Z_3)\M_{|}(H_{2a}/\Z_3)\tag\eak\\
\M(H_{2a}/\Z_6)&=\M_{|}(H_{2a}/\Z_6)^2\tag\eal\\
\endalign
$$
The strikingly uniform appearance of these identities (which as we have seen are equivalent to equations (\eaa)--(\ead)) invites one to search for an explanation as to why they hold. 

One way of trying to understand this is to place these identities in a larger context, by finding suitable generalizations of them. We discuss this in the next section.

We end this introduction with a brief account on how the group of identities (\eai)--(\eal), and the generalizations (\eba) and (\ebb) discussed in the next section, were discovered by the author. The factorization theorem of \cite{\FT} for perfect matchings of symmetric planar graphs readily implies that $\M_{-}(H_{a,a,2b})$ is a divisor of $\M(H_{a,a,2b})$, or, in the language of plane partitions, that $TC(a,a,2b)|P(a,a,2b)$. This raises the natural question of expressing their ratio in terms of related objects. A comparison of the resulting quantities for small values of $a$ and $b$ with the number of corresponding plane partitions in various symmetry classes quickly revealed that this ratio is equal to $S(a,a,2b)$, i.e. to $\M_{|}(H_{a,a,2b})$. This is how we found (\eai). Identity (\ead) was well-known, and its translation into the above context led to (\eal) (and was proved directly, without separately evaluating both sides, in \cite{\csts}). It then seemed natural to bring in various other rotational symmetries of the hexagon, and this way we were led to the uniform group of identities (\eaj)--(\eal). Identity (\eac) follows from Stembridge's results in \cite{\StemTS} by the factorization theorem of \cite{\FT}.

In the light of our previous work \cite{\ppone}, which was initially motivated by the desire to find a generalization of MacMahon's theorem on boxed plane partitions \cite{\MacM}, I then considered the regions described in the next section, and conjectured that the natural extensions of (\eai) and (\eaj) for them (namely, (\eba) and (\ebb)) also hold. This generalization of (\eai) was proved in \cite{\fakt}. The proof of the generalization of (\eaj) is the main result of the present paper.

%\pagebreak

\mysec{2. Statement of main result}

\topinsert
\centerline{\mypic{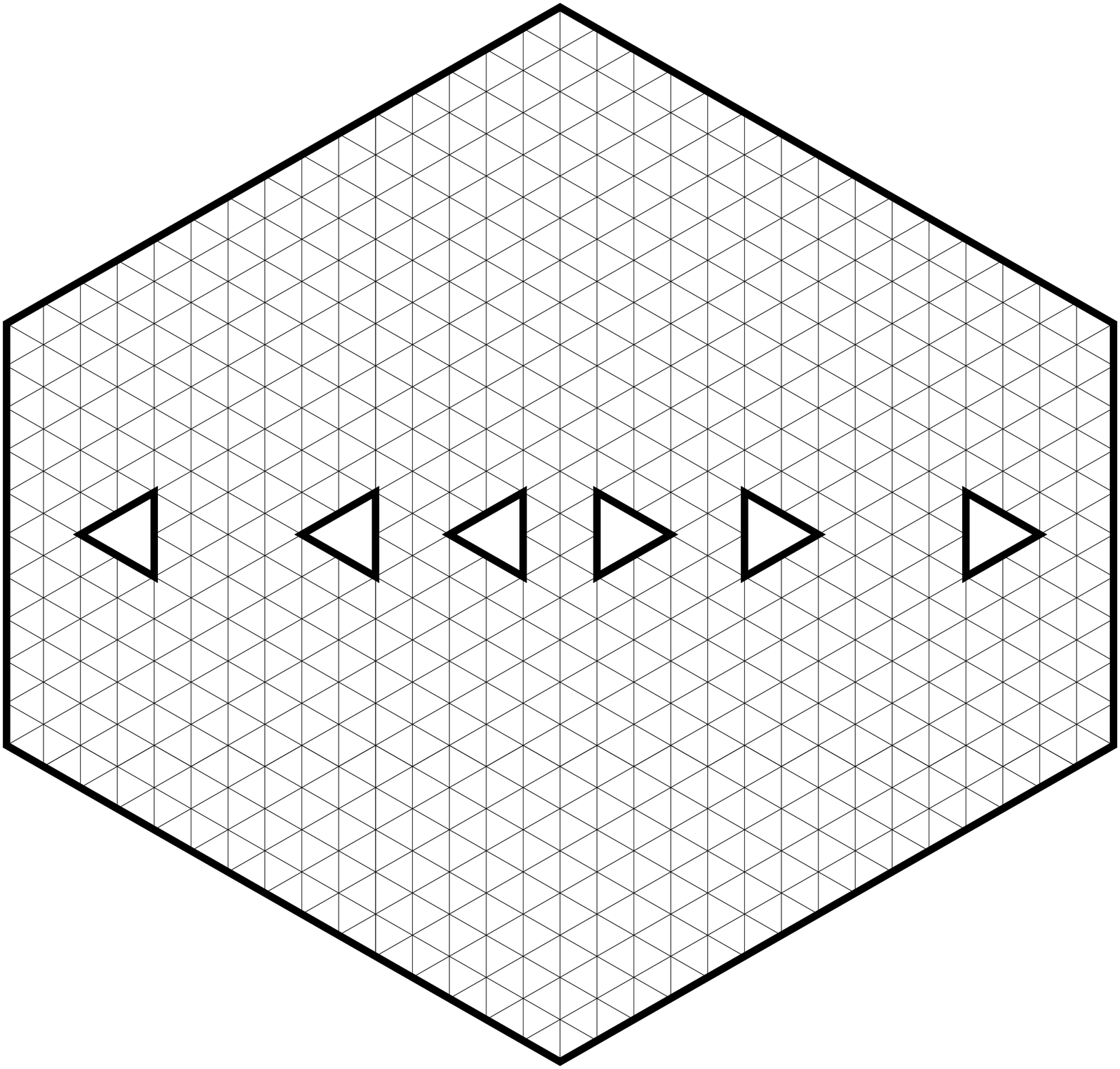}}
\centerline{{\smc Figure~\fba. {\rm The hexagon with holes $H_{15,10}(2,5,7)$.}}}
\endinsert

In \cite{\ppone} we extended MacMahon's theorem --- thought of as the enumeration of lozenge tilings of a hexagon on the triangular lattice --- by introducing, in the case when the hexagon has a symmetry axis, a certain collection of triangular holes straddling this axis, in such a way that the lozenge tilings of the resulting regions are still enumerated by simple product formulas. Inspired by that construction, it is natural to consider the family of regions $H_{a,2b}(k_1,\dotsc,k_s)$ defined as follows.

Let $a,b,s$ and $k_1,\dotsc,k_s$ be positive integers with $0< k_1<k_2<\dots<k_s\le a/2$.
Denote by $H_{a,2b}(k_1,k_2,\dots,k_s)$ the region obtained from the hexagon with side lengths $a,a,2b,a,a,2b$ (in clockwise order, starting from the northwestern side) by removing the following $2s$ triangles of side-length two from along its horizontal symmetry axis: $s$ left-pointing such triangles, with vertical sides at distances $2k_1,2k_2,\dots,2k_s$ from the left side of the hexagon (in units equal to $\sqrt{3}$ times the lattice spacing), and their reflections across the vertical symmetry axis of the hexagon.
%, out of which left-pointing triangles of side length
%two have been removed along the horizontal symmetry axis of the
%hexagon, their vertical sides being at distances
%$2k_1,2k_2,\dots,2k_l$ from the left side of the hexagon (in units equal to $\s%$rt{3}$ times the lattice spacing), and
%symmetrically corresponding right-pointing triangles of side length
%two have been removed along the horizontal symmetry axis, 
%their vertical sides being at lattice distances
%$2k_1,2k_2,\dots,2k_l$ from the right side of the hexagon.
Figure~\fba\ shows the region $H_{15,10}(2,5,7)$.

Consider also the following variant of the above regions, obtained by removing from their center a horizontal lattice rhombus of odd side-length (as explained in \cite{\fakt}, removing a central rhombus of even side-length does not lead to new regions, but to ones that are equivalent to certain special cases of the $H_{a,2b}(k_1,k_2,\dots,k_s)$'s).
 
Let $a,b,s$ and $k_1,\dotsc,k_s$ be positive integers with $0< k_1<k_2<\dots<k_s\le a/2$.
One readily sees that the region
$H_{a,2b}(k_1,k_2,\dotsc,k_s)$ has a horizontal lattice rhombus of
odd side-length at its center precisely if $a$ is odd. 
%Assume that this is so, and let 
For any positive integer $x\leq a$, denote therefore by $H_{2a-1,2b}(k_1,k_2,\dotsc,k_s;2x-1)$ the region obtained from $H_{2a-1,2b}(k_1,k_2,\dotsc,k_s)$ by removing from its center the horizontal lattice rhombus of side $2x-1$ (see Figure~{\fbb} for an example).

\topinsert
\centerline{\mypic{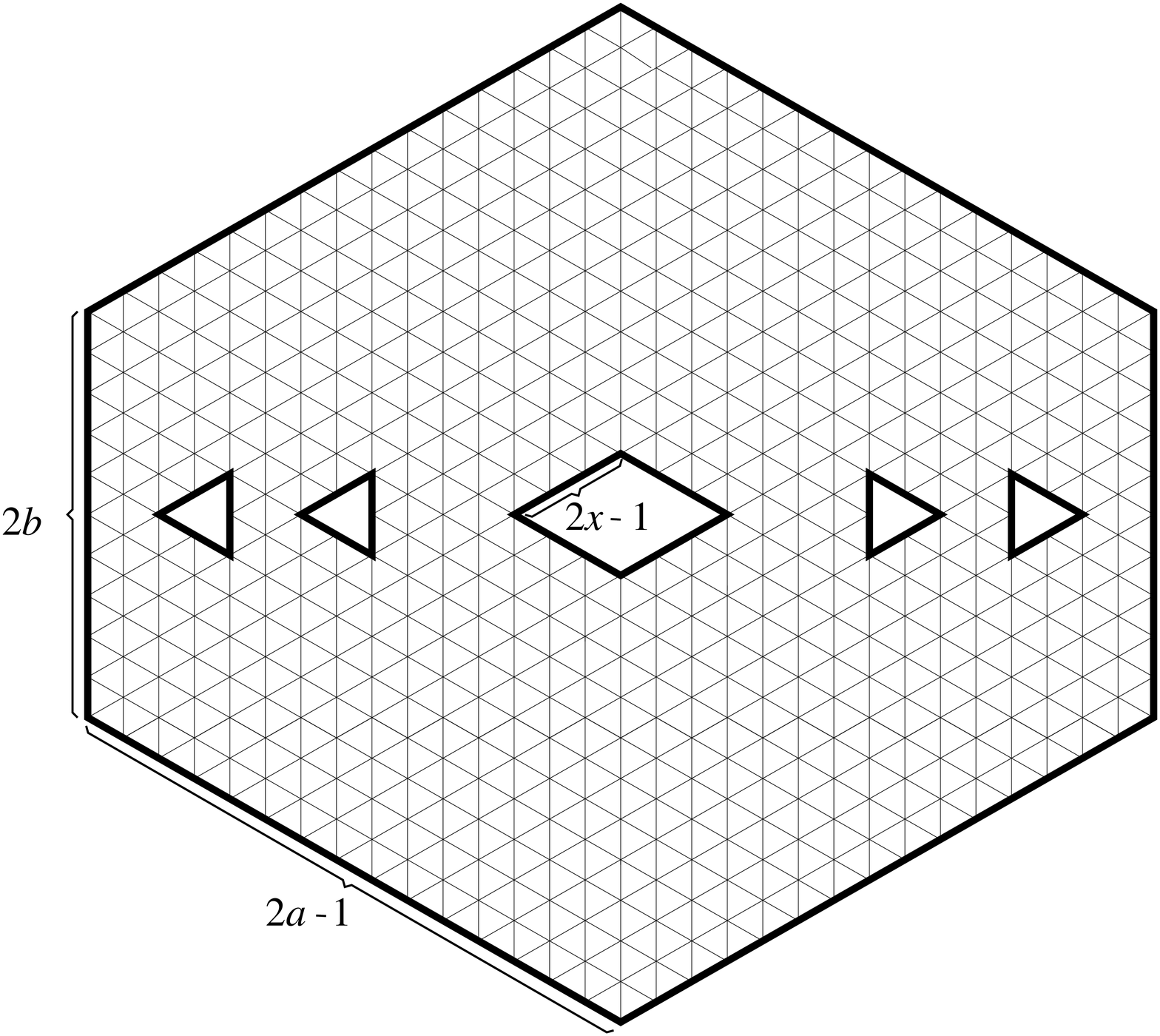}}
\centerline{{\smc Figure~\fbb. {\rm The region $H_{15,10}(2,4;3)$.}}}
\endinsert

As proved in \cite{\fakt}, it turns out that for all positive integers $a,b,s,x$ with $x\leq a$ and all integers $k_1,k_2,\dotsc,k_s$ with $0< k_1<k_2<\dots<k_s\le a/2$, the following natural extensions of (\eai) to the above two families of regions hold:
$$
\M(H_{a,2b}(k_1,k_2,\dots,k_s))
=
\M_{-}(H_{a,2b}(k_1,k_2,\dots,k_s))
\M_{|}(H_{a,2b}(k_1,k_2,\dots,k_s))
\tag\eba
$$
and
$$
\spreadlines{3\jot}
\align
\!\!\!\!\!\!\!\!\!\!\!\!\!\!\!\!\!\!\!\!\!\!\!\!
&
\M(H_{2a-1,2b}(k_1,k_2,\dots,k_s;2x-1))
=
\M_{-}(H_{2a-1,2b}(k_1,k_2,\dots,k_s;2x-1))
\\
&\ \ \ \ \ \ \ \ \ \ \ \ \ \ \ \ \ \ \ \ \ \ \ \ \ \ \ \ \ \ \ \ \ \ \ \ \ \ \ \ \ \ \ \ 
\times
\M_{|}(H_{2a-1,2b}(k_1,k_2,\dots,k_s;2x-1)).
\tag\ebb
\endalign
$$
Note that all quantities involved in (\eaj) are also defined for the generalizations $H_{a,2b}(k_1,k_2,\dots,k_s)$ and $H_{2a-1,2b}(k_1,k_2,\dots,k_s;2x-1)$ of the hexagonal regions. The main result of this paper is to prove that for the generalized regions, these quantities are still related by the same equations. More precisely, the following result holds.

\proclaim{Theorem \tba} For all positive integers $a,b,s,x$ with $x\leq a$ and all integers $k_1,k_2,\dotsc,k_s$ with $0< k_1<k_2<\dots<k_s\le a/2$, we have
%$$
%\M(H_{a,2b}(k_1,k_2,\dotsc,k_l)/\Z_2)
%=
%\M_{|}(H_{a,2b}(k_1,k_2,\dotsc,k_l)/\Z_2)^2.
%\tag\ebb
%$$
$$
\M_{\odot}(H_{a,2b}(k_1,k_2,\dotsc,k_s))
=
\M_{\odot,|}(H_{a,2b}(k_1,k_2,\dotsc,k_s))^2,
\tag\ebc
$$
and
$$
\M_{\odot}(H_{2a-1,2b}(k_1,k_2,\dotsc,k_s;2x-1))
=
\M_{\odot,|}(H_{2a-1,2b}(k_1,k_2,\dotsc,k_s;2x-1))^2,
\tag\ebd
$$
where $\odot$ denotes symmetry with respect to the center of the region (equivalently, invariance under rotation by $180^\circ$ around its center).

\endproclaim

%\pagebreak

\mysec{3. Proof of Theorem {\tba}}

We begin with the observation that we may assume without loss of generality that $k_1\neq1$. Indeed, if $k_1=1$, the holes farthest from the center touch the vertical sides of $H_{a,2b}(k_1,\dotsc,k_s)$, causing two rows of forced lozenges along the left and right boundaries. Upon removing these forced lozenges, the leftover region is a smaller region of type $H_{a,2b}(k_1,\dotsc,k_s)$, in which the value of $k_1$ is not equal to 1. Thus the tilings of the original region are naturally identified with the tilings of a smaller region of the same type, but with $k_1\neq1$. Furthermore, under this identification, vertically symmetric tilings of the original region are mapped to vertically symmetric tilings of the smaller region, and similarly for centrally symmetric and vertically-symmetric-and-centrally-symmetric tilings. This implies the statement at the beginning of this paragraph.

%the case $k_1=1$ of Theorem {\tba} follows from the case $k_1>1$. 

We prove (\ebc) and (\ebd) by separately evaluating their left and right hand sides, and verifying that they agree.
In proving (\ebc), the details of the arguments are slightly different depending on the parity of $a$.

\topinsert
\centerline{\mypic{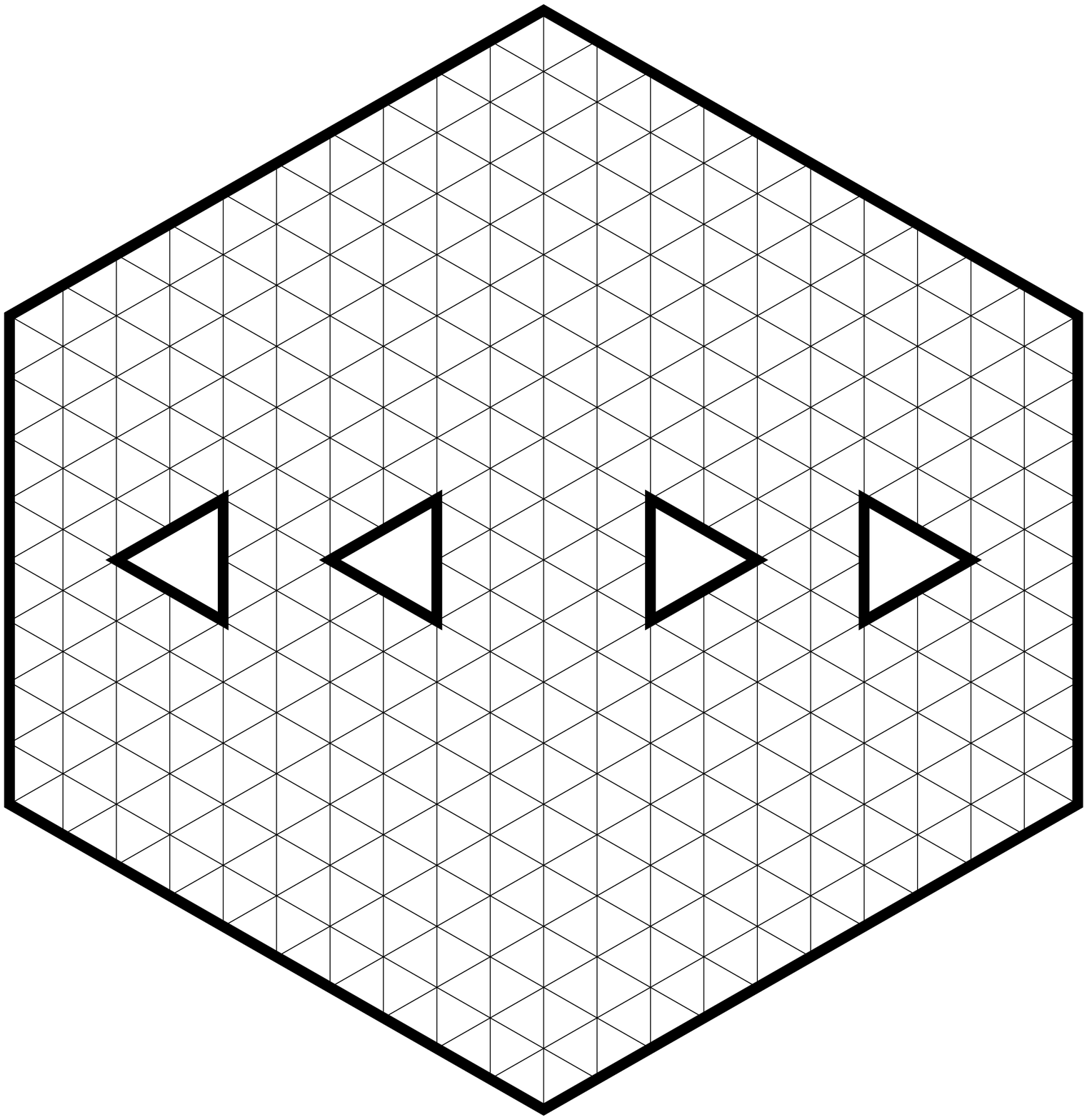}}
\centerline{{\smc Figure~\fca. {\rm The region $H_{10,8}(2,4)$.}}}
\endinsert

\topinsert
\twoline{\mypic{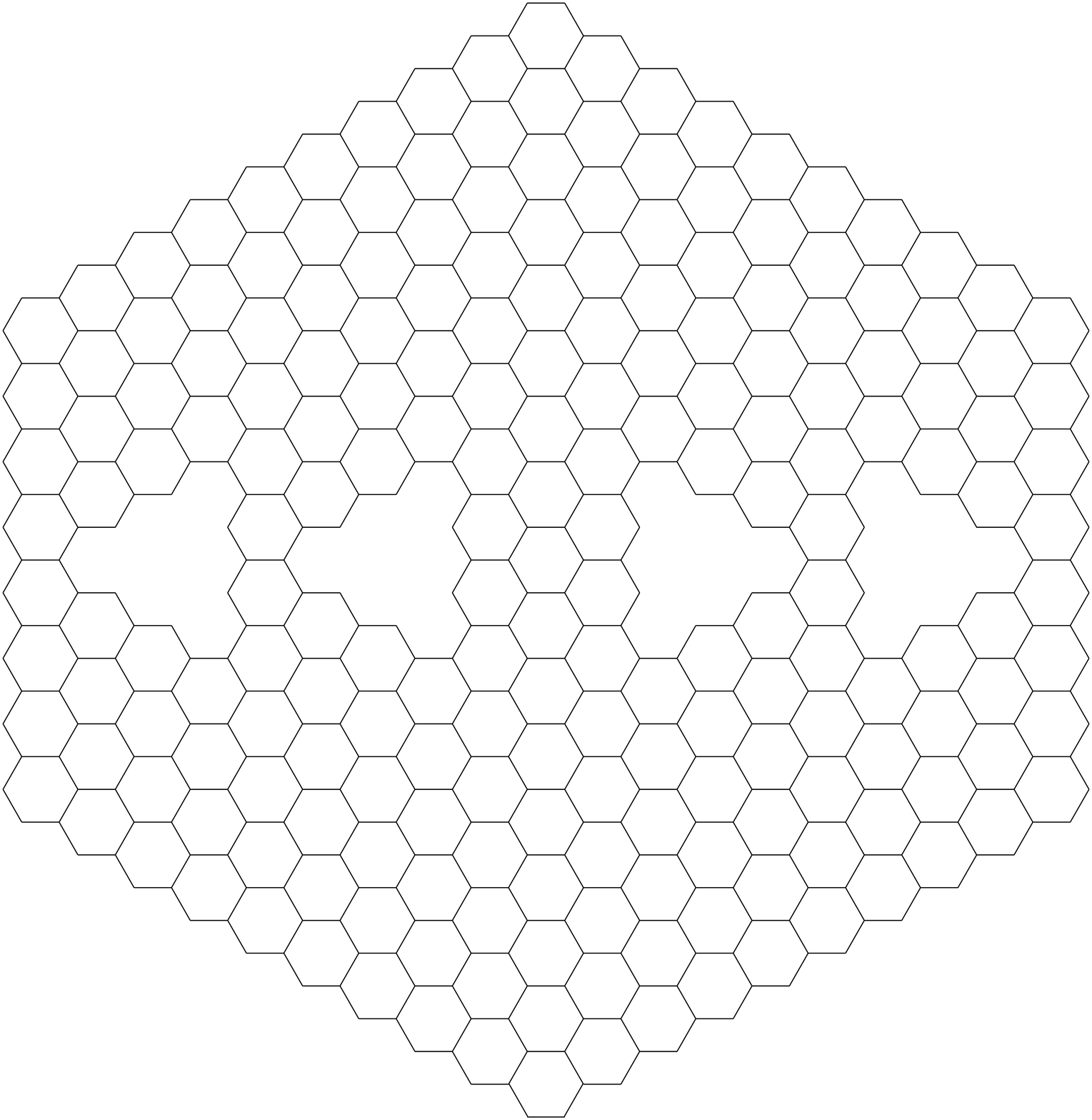}}{\mypic{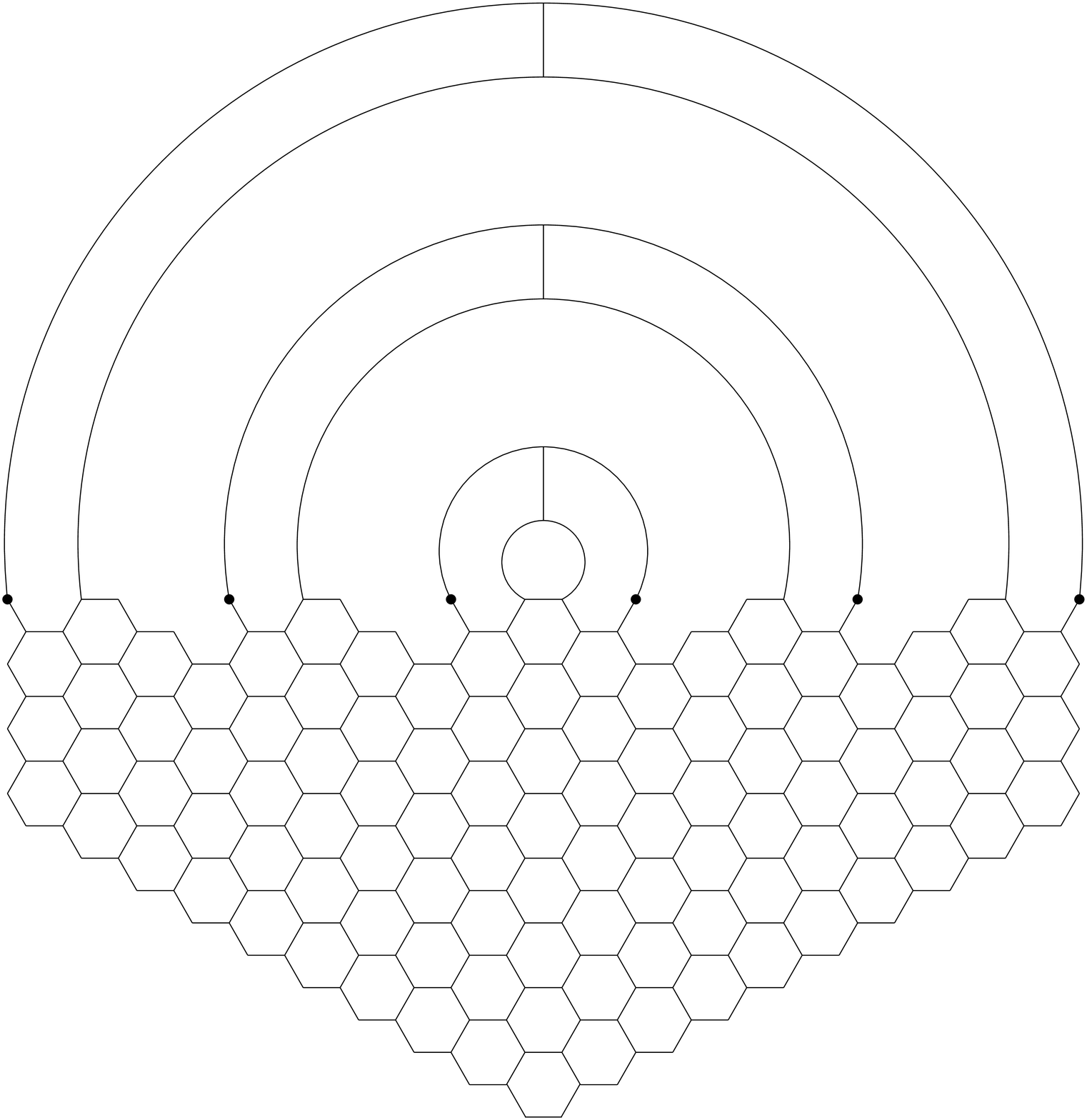}}
\medskip
\centerline{{\smc Figure~{\fcb}. {\rm  The dual graph of $H_{10,8}(2,4)$, and its quotient by $180^\circ$ rotation.}}}
\endinsert

{\it Case 1: $a$ is even} (an illustrative example is shown in Figure {\fca}). For clarity and notational simplicity, throughout this case we write $2a$ instead of $a$. 

As a consequence of its definition, the left hand side of (\ebc) is equal to the number of perfect matchings of the quotient of the dual graph of $H_{2a,2b}(k_1,k_2,\dotsc,k_s)$ under rotation by $180^\circ$ around its center (see Figure {\fcb} for an illustration of the dual graph and its quotient). One readily sees that this quotient graph can be embedded symmetrically in the plane (this is illustrated in the picture on the right in Figure {\fcb}; the black dots indicate vertices of degree two that may not otherwise be apparent). 
% (see Figure {\fba} for an example --- the region, its dual and the quotient).
It can be easily checked that the variant of the factorization theorem \cite{\FT,Theorem\,1.2} described in \cite{\FT, Proof\ of\ Theorem\,7.1} can be applied to this graph. 
One obtains that the number of matchings of
(the dual graph of) $H_{2a,2b}(k_1,k_2,\dotsc,k_s)$ equals $2^{a-s}$ times the matching generating function\footnote{ The matching generating
function of a graph is the sum of the weights of all its perfect matchings, where the
weight of a matching is the product of the weights of its edges.} of the subgraph $K_{2a,2b}(k_1,k_2,\dotsc,k_s)$ (illustrated on the left in Figure {\fcc}) 
obtained by deleting its top $2a-2s$ edges immediately to the right of the symmetry axis, 
and changing the weight of the $a-s$ edges along the symmetry axis to 1/2 (the resulting graph has been redrawn in Figure {\fcc} so that it is a subgraph of the hexagonal lattice; the edges weighted by $1/2$ are marked).

\topinsert
\twoline{\mypic{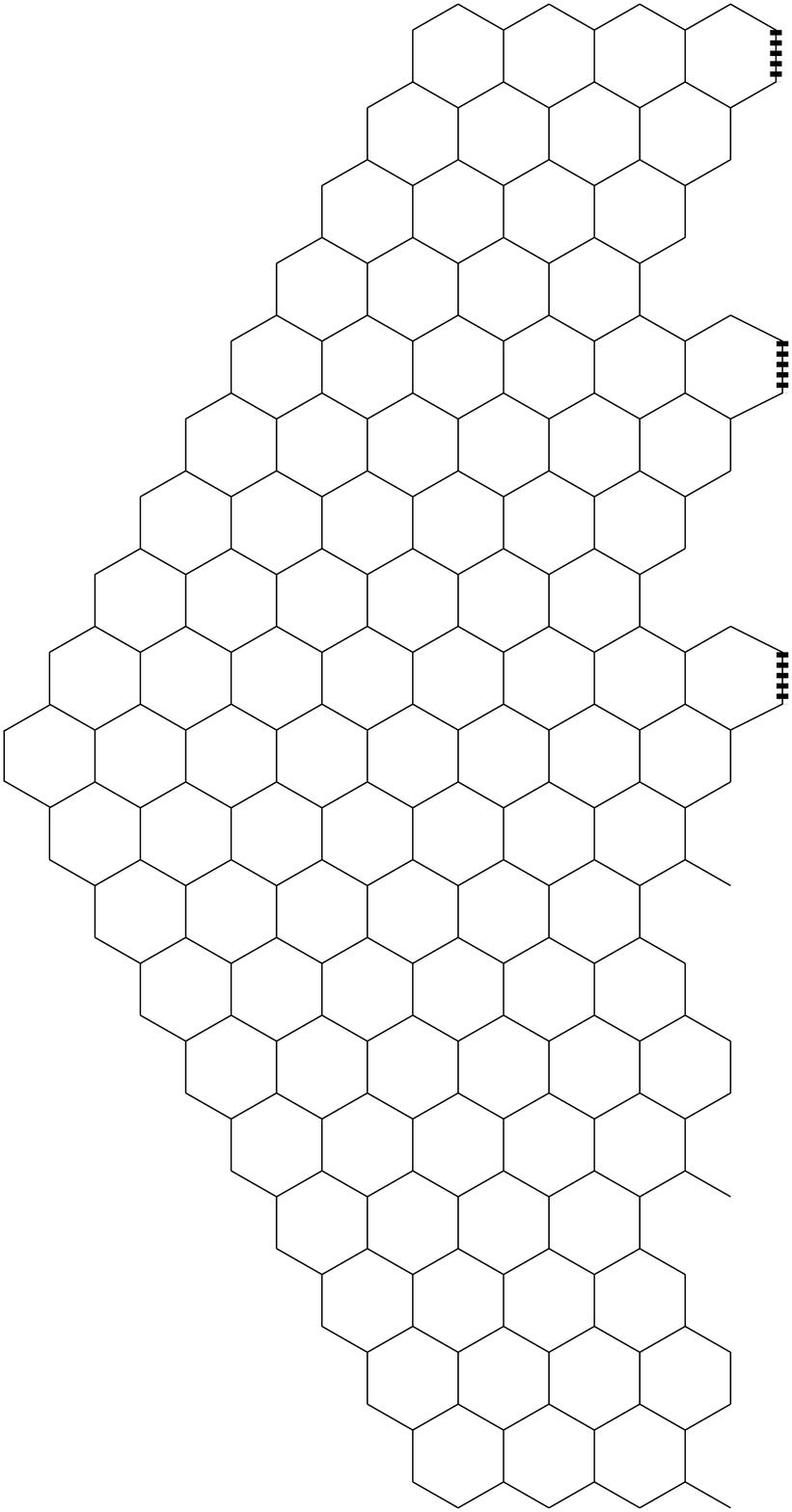}}{\mypic{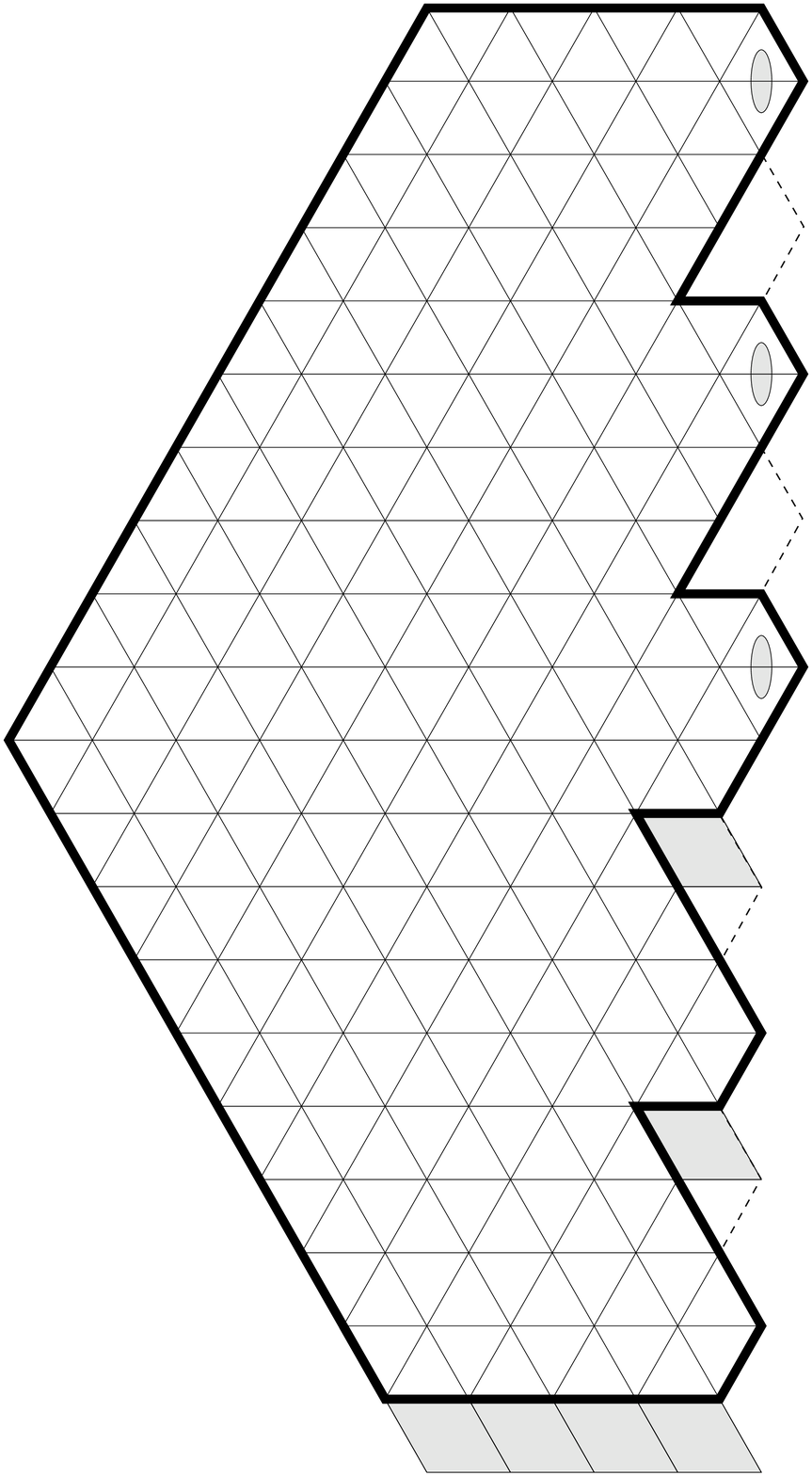}}
\medskip
\centerline{{\smc Figure~{\fcc}. {\rm  The graph $K_{10,8}(2,4)$ resulting by applying the }}}
\centerline{{{\rm factorization theorem to $H_{10,8}(2,4)$, and its dual region.}}}
\endinsert

Using  
again the duality between matchings and lozenge tilings, the matchings of the graph $K_{2a,2b}(k_1,k_2,\dotsc,k_s)$ can be
identified with tilings of its dual region --- which we also denote by $K_{2a,2b}(k_1,k_2,\dotsc,k_s)$, according to our notation-simplifying convention ---  illustrated on the right in Figure {\fcc}. Consider the $a-s$ easternmost tile positions in this region
along its northeastern boundary (they are indicated by shaded ellipses in Figure {\fcc}). In a
tiling of this region, weight each tile occupying one of these positions by 1/2, and all
others by~1.
%; write $L^*(R_n)$ be the tiling generating function of $R_n$ under this weighting. 
Then the bijection between matchings of $K_{2a,2b}(k_1,k_2,\dotsc,k_s)$ and
tilings of its dual region is weight-preserving. Therefore, one obtains
$$
\M(H_{2a,2b}(k_1,k_2,\dotsc,k_s))=2^{a-s}\M(K_{2a,2b}(k_1,k_2,\dotsc,k_s)).
\tag\eca
$$
However, the resulting region $K_{2a,2b}(k_1,k_2,\dotsc,k_s)$ belongs to the family of $\bar{R}$-regions defined in \cite{\ppone,Section\,2}, whose tilings are enumerated by \cite{\ppone,Proposition\,2.1}. Namely, one readily verifies that $K_{2a,2b}(k_1,k_2,\dotsc,k_s)$ is precisely the region\footnote{ We denote by $[n]$ the list $[1,2,\dotsc,n]$.}
$$
\bar{R}_{\left[a-1\right]\setminus\left[a-k_1,\dotsc,a-k_s\right],
\left[a\right]\setminus\left[a-k_1+1,\dotsc,a-k_s+1\right]}(b)
\tag\ecb
$$
defined in \cite{\ppone, Section\,2} (the region on the right in Figure {\fcc} is $\bar{R}_{[2,4],[1,3,5]}(4)$ in the language of \cite{\ppone}, as on its right it contains the 2nd and 4th ``bumps'' of the lower zig-zag line supporting its boundary, the 1st, 3rd and 5th bumps of the upper zig-zag, and the length of the base is 4).

Define the lists ${\bold l}=[l_1,\dotsc,l_{a-s-1}]$ and  ${\bold q}=[q_1,\dotsc,q_{a-s}]$ (with $l_1<\cdots<l_{a-s-1}$ and $q_1<\cdots<q_{a-s}$) by
$$
{\bold l}:=\left[a-1\right]\setminus\left[a-k_1,\dotsc,a-k_s\right]
\tag\ecc
$$
and
$$
{\bold q}:=\left[a\right]\setminus\left[a-k_1+1,\dotsc,a-k_s+1\right].
\tag\ecd
$$
Then, using \cite{\ppone,Proposition\,2.1}, we obtain from (\eca) that
$$
\spreadlines{3\jot}
\align
&
\M(H_{2a,2b}(k_1,k_2,\dotsc,k_s))=\frac12
\prod_{i=1}^{a-1}\frac{1}{(2l_i-1)!}\prod_{i=1}^{a}\frac{1}{(2q_i)!}
\\
&\ \ \ \ \ \ \ 
\times
\frac{\prod_{1\leq i<j\leq a-s-1}(l_j-l_i)\prod_{1\leq i<j\leq a-s}(q_j-q_i)}
{\prod_{i=1}^{a-s-1}\prod_{j=1}^{a-s}(l_i+q_j)}
\,Q_{{\bold l},{\bold q}}(a+b-s),
\tag\ece
\endalign
$$
where the polynomial $Q_{{\bold l},{\bold q}}(x)$ is given by\footnote{ In the first line of (\ecf), the bases are incremented by 1 from each factor to the next, while the exponents are incremented by 1 from each factor to the next until the middle, and then decremented by one at each step.}
$$
\spreadlines{3\jot}
\align
&
Q_{{\bold l},{\bold q}}(x):=
\left((x+1) (x+2)^2 \cdots (x+2a-2s-2)^2 (x+2a-2s-1)\right)^2
\\
&\ \ \ \ \ \ \ \ \ \ \ \ \ \ \ 
\times
\prod_{i=1}^{a-s-1}\prod_{j=1}^{l_i-i}(x-i-j+a-s+1)(x+i+j+a-s-1)
\\
&\ \ \ \ \ \ \ \ \ \ \ \ \ \ \ 
\times
\prod_{i=1}^{a-s}\prod_{j=1}^{q_i-i}(x-i-j+a-s)(x+i+j+a-s).
\tag\ecf
\endalign
$$

\topinsert
\centerline{\mypic{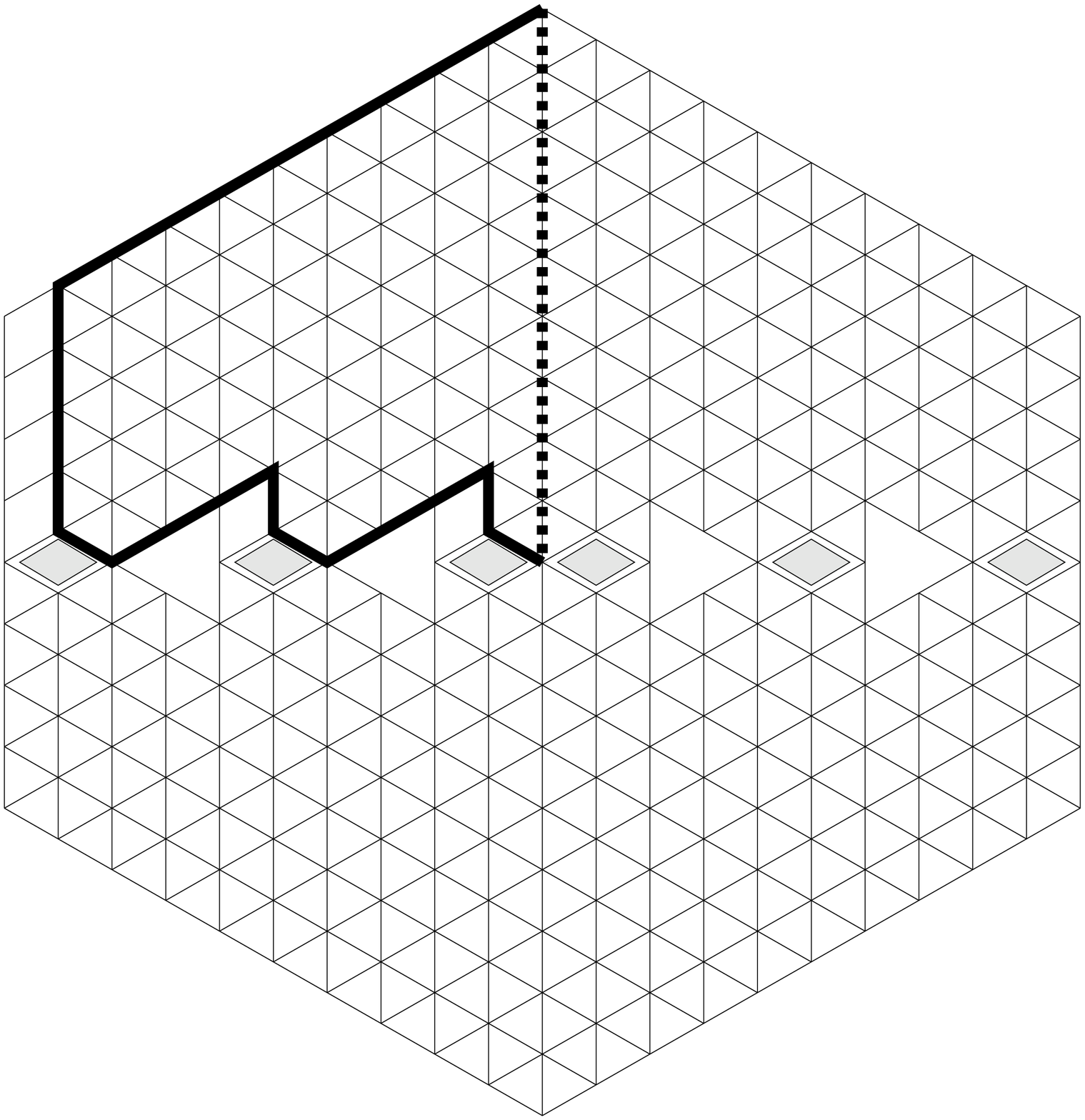}}
\centerline{{\smc Figure~\fcd. {\rm The region $D_{5,4,-1}^{1,3,5}$.}}}
\endinsert

We now turn to the evaluation of the right hand side of (\ebc). Note that for a lozenge tiling of $H_{2a,2b}(k_1,\dotsc,k_s)$, to be invariant under both rotation by $180^\circ$ and reflection across the vertical is equivalent to being invariant under reflection across both the horizontal and vertical symmetry axes. Denote these symmetry axes by $\ell_h$ and $\ell_v$, respectively. 

Any lozenge tiling of $H_{2a,2b}(k_1,\dotsc,k_s)$ which is symmetric across the horizontal must contain the $2a-2s$ lozenges indicated by a shading in Figure {\fcd}. When removing these lozenges, $H_{2a,2b}(k_1,\dotsc,k_s)$ gets disconnected into two congruent regions, one above and one below $\ell_h$; the horizontally symmetric tilings of $H_{2a,2b}(k_1,\dotsc,k_s)$ are in bijection with the tilings of say the upper region.

Then the horizontally {\it and vertically} symmetric tilings of $H_{2a,2b}(k_1,\dotsc,k_s)$ are in bijection with the tilings of the upper region which are symmetric across $\ell_v$. In turn, after removing the forced lozenges from the upper region, these tilings are readily seen to be in bijection with the lozenge tilings of the left half of the upper region, provided its boundary along $\ell_v$ is considered free --- i.e., in its tilings lozenges can protrude out halfway across it (this region is indicated in Figure {\fcd} by the thick contour; the free portion of its boundary is indicated by the dashed line). However, under our assumption that $k_1>1$, if we set
$$
\{i_1,\dotsc,i_{a-s}\}:=\{1,\dotsc,a\}\setminus\{a-k_1+1,\dotsc,a-k_s+1\},
$$
with $i_1<\cdots<i_{a-s}$, this is precisely the region $D_{a,b,-1}^{i_1,\dotsc,i_{a-s}}$ defined in \cite{\ranglep,Section\,3}. It follows then from the above and \cite{\ranglep,Proposition\,3.1} that we have
$$
\spreadlines{3\jot}
\align
\M_{\odot,|}(H_{2a,2b}(k_1,\dotsc,k_s))
&=\M_f(D_{a,b,-1}^{i_1,\dotsc,i_{a-s}}))\\
&=
\prod_{j=1}^{a-s}{a+b+i_j-1 \choose 2i_j-1}
\prod_{1\leq j<k\leq a-s}\frac{i_k-i_j}{i_j+i_k-1}.
\tag\ecg
\endalign
$$
It is routine to verify that the expression given by (\ece)--(\ecf) is precisely the square of the expression~(\ecg). This proves (\ebc) for even $a$.

\topinsert
\centerline{\mypic{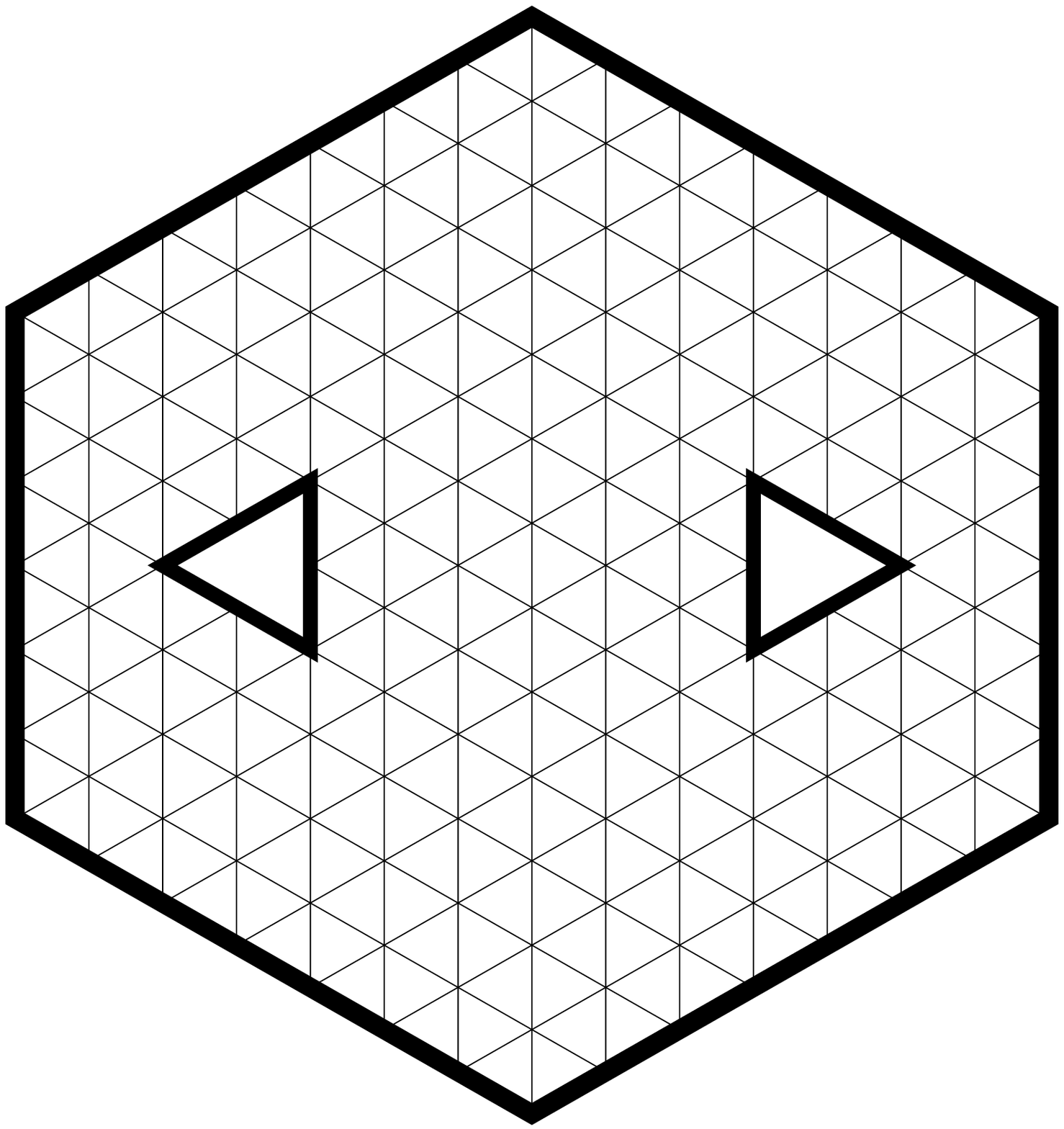}}
\centerline{{\smc Figure~\fce. {\rm The region $H_{7,6}(2)$.}}}
\endinsert

\topinsert
\twoline{\mypic{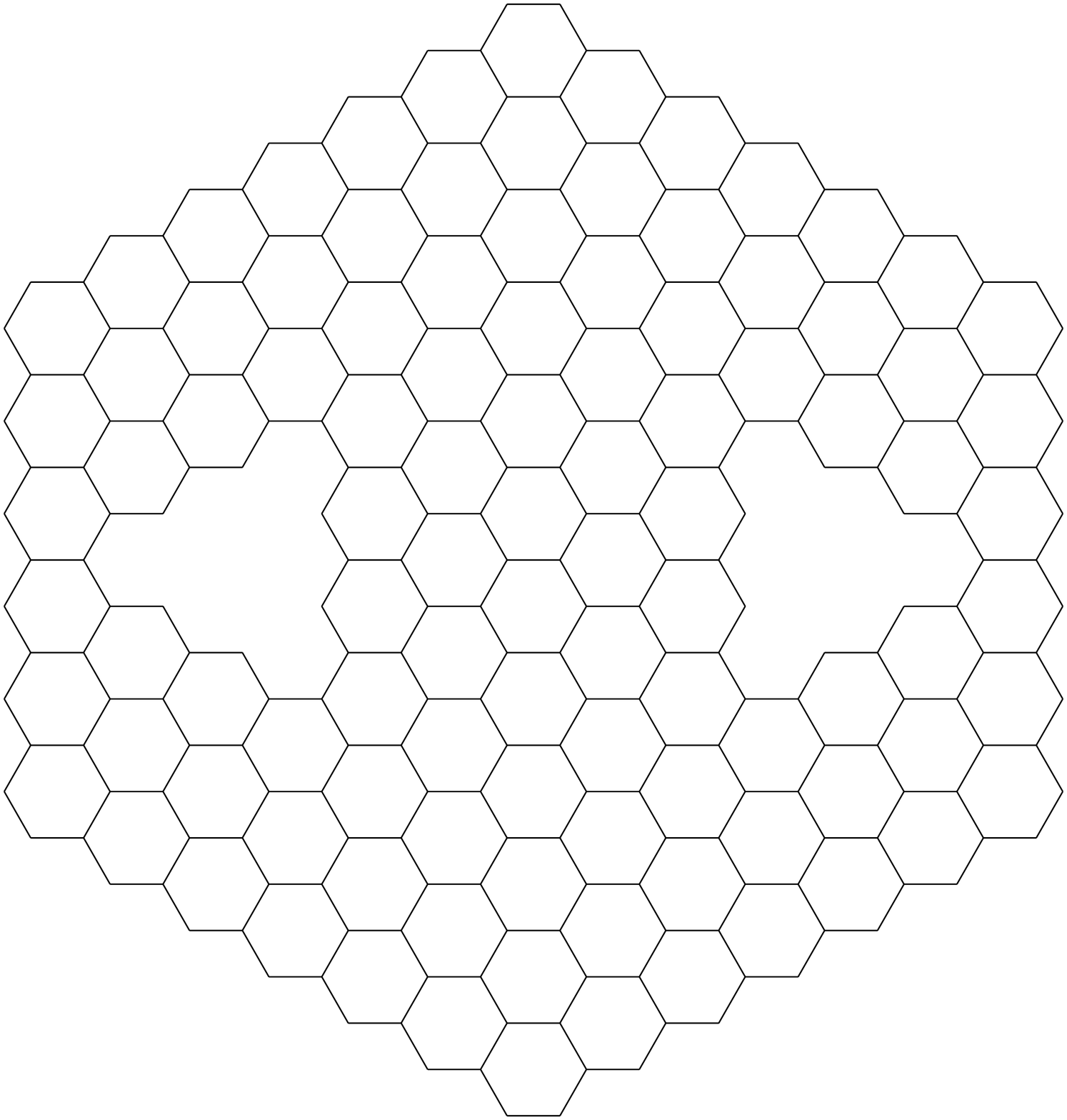}}{\mypic{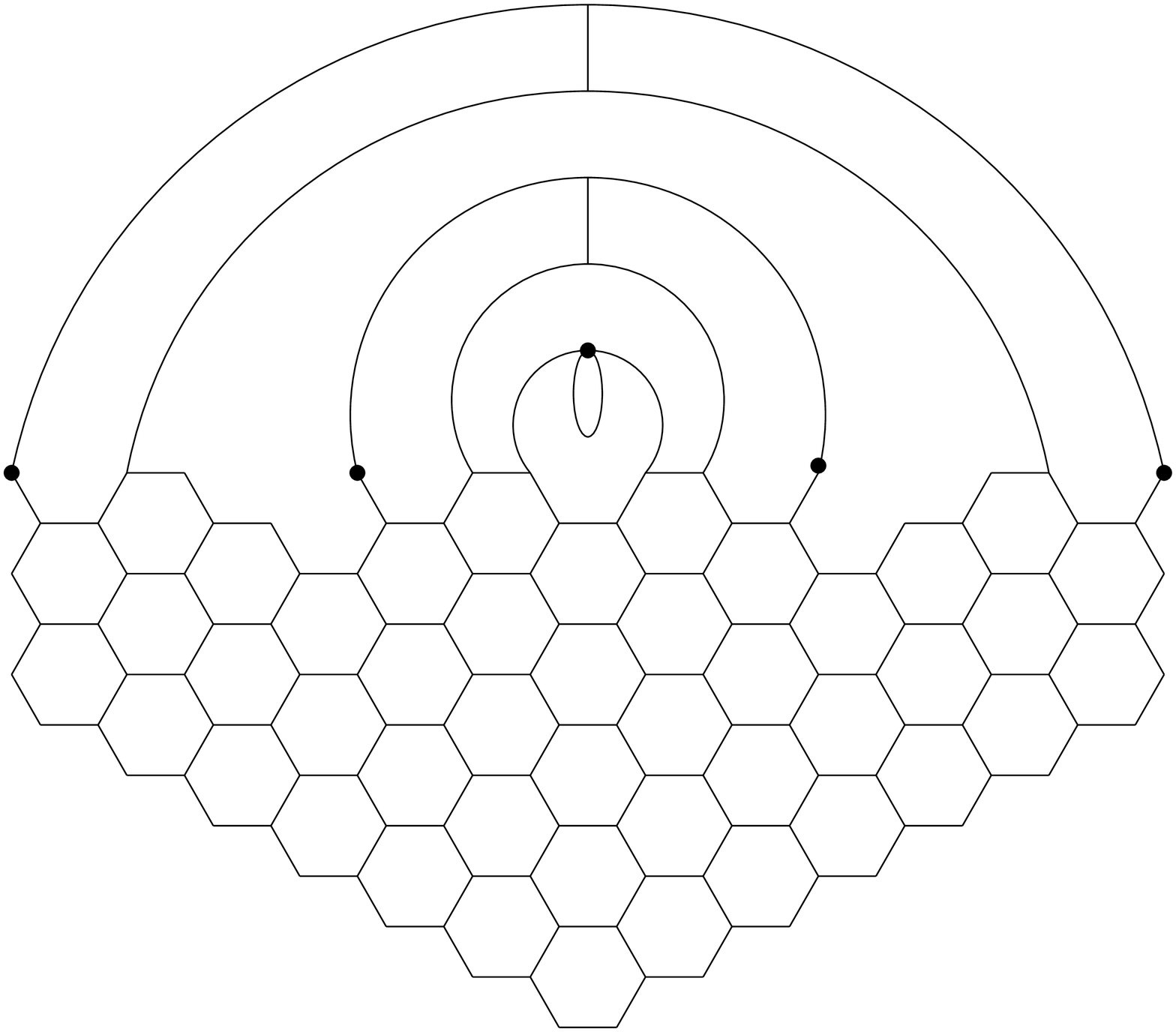}}
\medskip
\centerline{{\smc Figure~{\fcf}. {\rm  The dual graph of $H_{7,6}(2)$, and its quotient by $180^\circ$ rotation.}}}
\endinsert

%For odd $a$ --- point out only the differences: mainly that we have to discard a loop from the quotient graph.

{\it Case 2. $a$ is odd} (an illustrative example is shown in Figure {\fce}). In analogy to the previous case, in the interest of clarity and notational simplicity, throughout this case we write $2a+1$ instead of $a$.

\topinsert
\twoline{\mypic{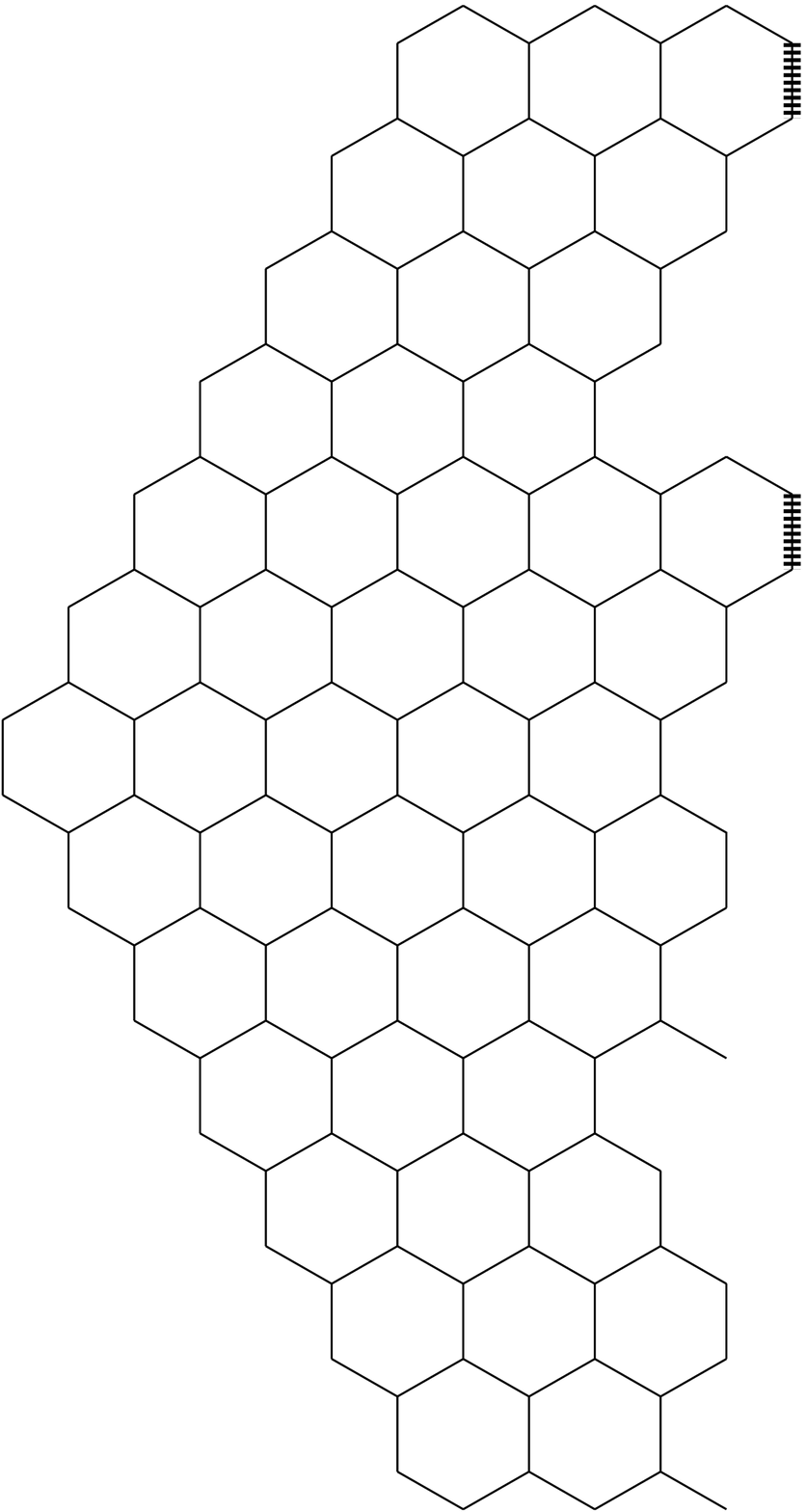}}{\mypic{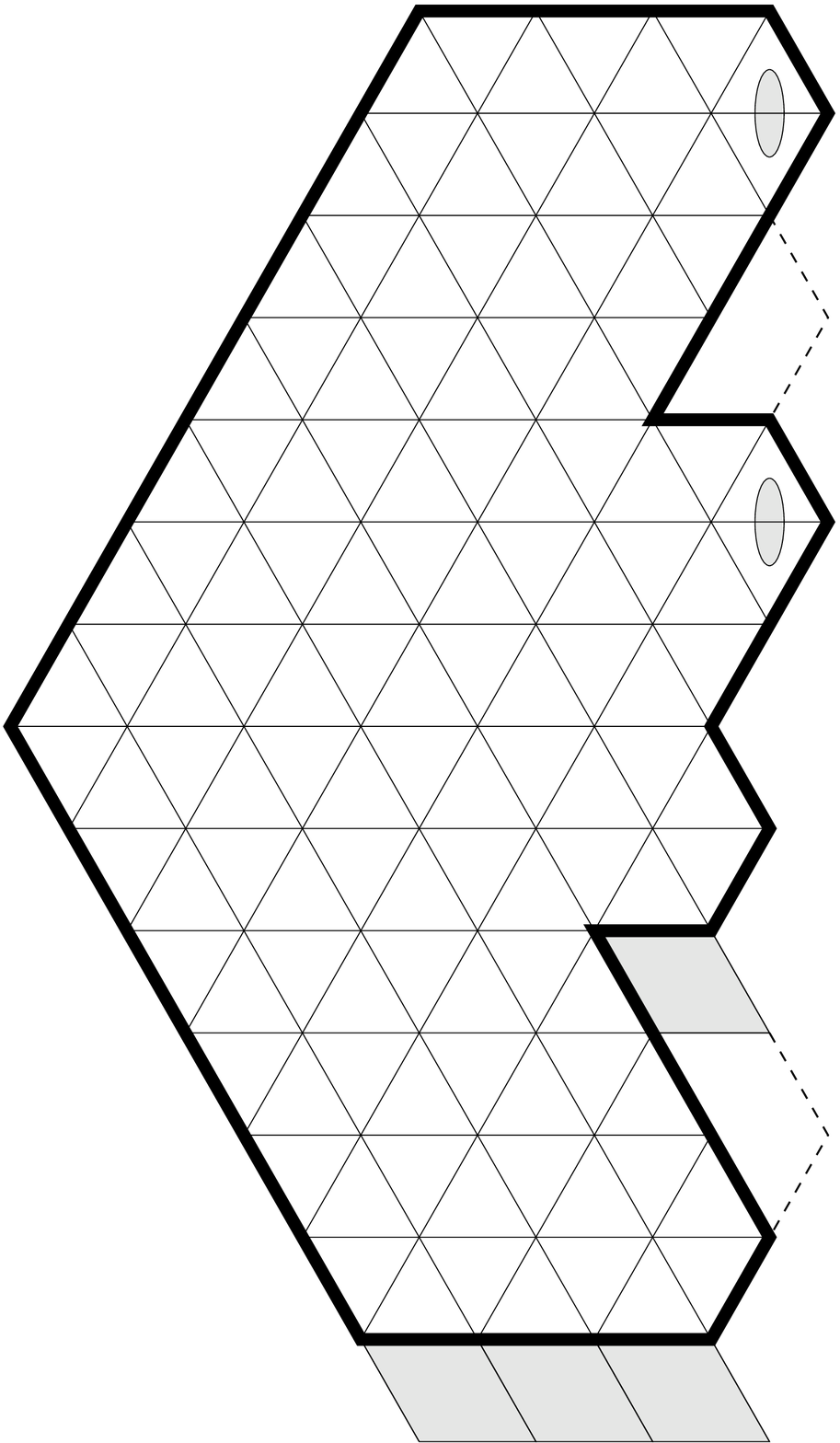}}
\medskip
\centerline{{\smc Figure~{\fcg}. {\rm  The graph $K_{7,6}(2)$ resulting by applying the }}}
\centerline{{{\rm factorization theorem to $H_{7,6}(2)$, and its dual region.}}}
\endinsert

The same reasoning as in Case 1 shows that 
$$
\M(H_{2a+1,2b}(k_1,k_2,\dotsc,k_s))=2^{a-s}\M(K_{2a+1,2b}(k_1,k_2,\dotsc,k_s)),
\tag\ech
$$
where $K_{2a+1,2b}(k_1,k_2,\dotsc,k_s)$ is the quotient graph of $H_{2a+1,2b}(k_1,k_2,\dotsc,k_s)$ with respect to rotation by $180^\circ$ around its center.

The main difference compared to Case 1 is that now this quotient graph contains a loop (see Figure {\fcf}). Since the quotient graph has an odd number of vertices, this loop must be part of all its perfect matchings. After removing from it the vertex at which this loop is based, applying the factorization theorem of \cite{\FT} to the resulting graph, and recognizing the resulting subgraph as being dual to a member of the family of $\bar{R}$-regions defined in \cite{\ppone,Section\,2} (see Figure {\fcg}), we are led to the equality
$$
\M(H_{2a+1,2b}(k_1,\dotsc,k_s))=2^{a-s}\M(\bar{R}_{[a]\setminus[a-k_1+1,\dotsc,a-k_s+1],[a]\setminus[a-k_1+1,\dotsc,a-k_s+1]}(b)).
\tag\eci
$$
Using the product expression that \cite{\ppone,Proposition\,2.1} supplies for the region on the right hand side above, we obtain that, with ${\bold q}$ given by (\ecd), we have  
$$
%\spreadlines{3\jot}
%\align
%&
\M(H_{2a+1,2b}(k_1,k_2,\dotsc,k_s))=
\prod_{i=1}^{a}\frac{1}{(2q_i-1)!(2q_i)!}
%\\
%&\ \ \ \ \ \ \ 
%\times
\frac{\left(\prod_{1\leq i<j\leq a-s}(q_j-q_i)\right)^2}
{\prod_{i=1}^{a-s}\prod_{j=1}^{a-s}(q_i+q_j)}
\,S_{{\bold q}}(a+b-s),
\tag\ecj
%\endalign
$$
where the polynomial $S_{{\bold q}}(x)$ is given by\footnote{ The pattern for the factors in the first line of (\eck) is the same as the one for the first line of (\ecf) --- which is explained in the previous footnote --- with the exception that the maximum value of the exponent is attained twice, for the middle two factors.} 
$$
\spreadlines{3\jot}
\align
&
S_{{\bold q}}(x):=
\left((x+1) (x+2)^2 \cdots (x+2a-2s-1)^2 (x+2a-2s)\right)^2
\\
&\ \ \ \ \ \ \ \ \ \ \ \ \ \ \ 
\times
\left(\prod_{i=1}^{a-s}\prod_{j=1}^{q_i-i}(x-i-j+a-s)(x+i+j+a-s)\right)^2.
\tag\eck
\endalign
$$

\topinsert
\twoline{\mypic{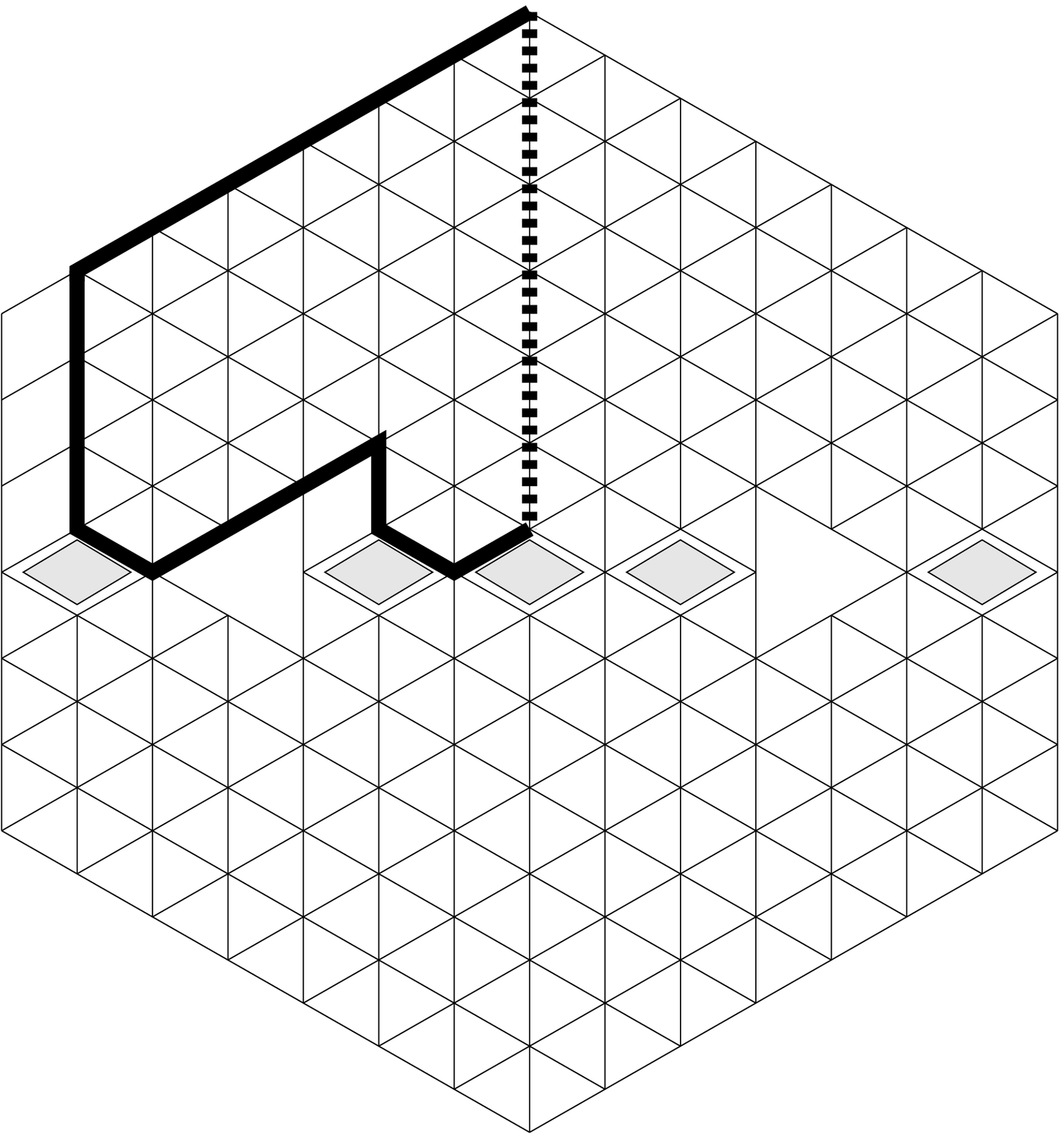}}{\mypic{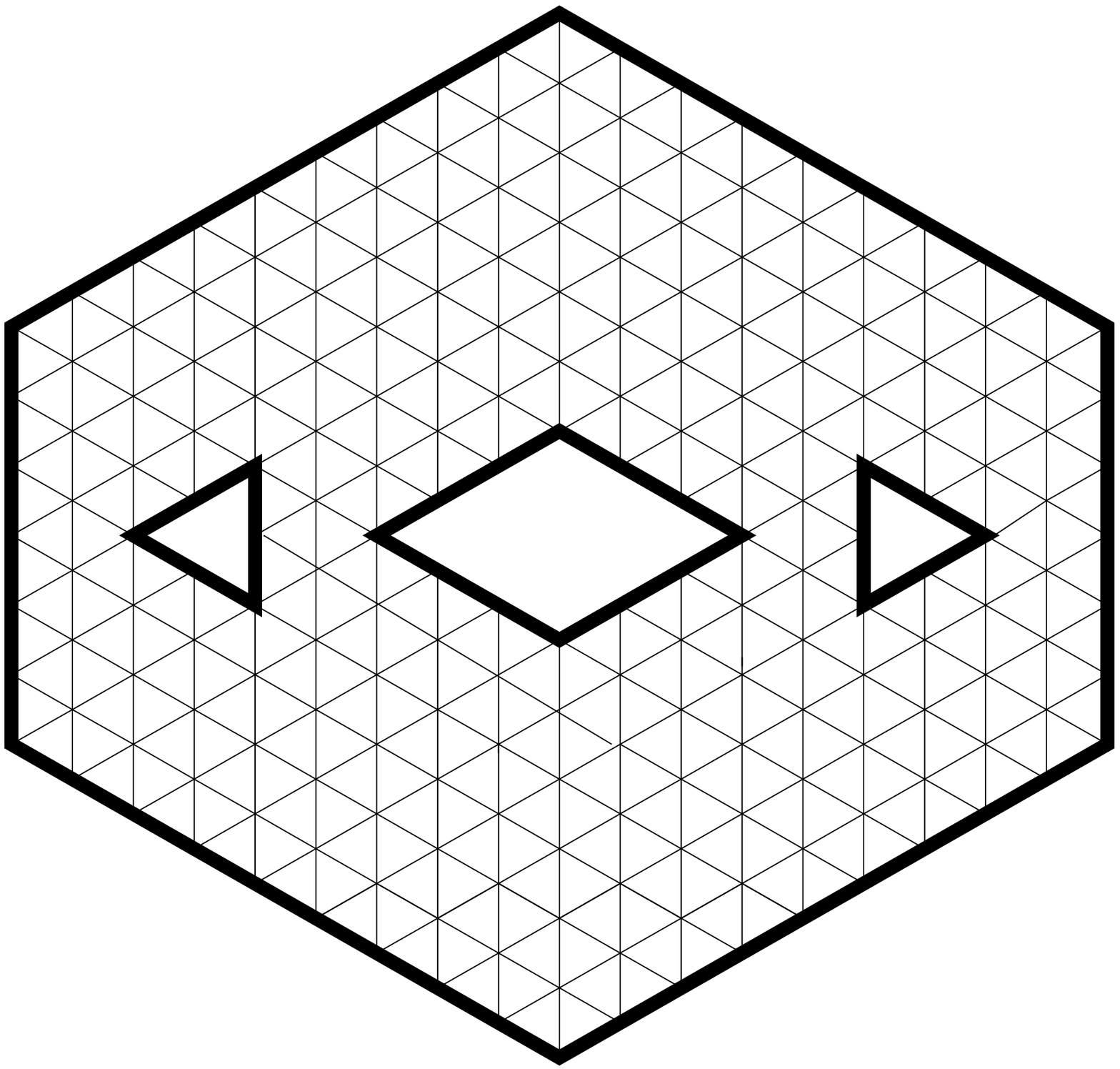}}
\twoline{{\smc Figure~\fch. {\rm The region $D_{3,3,0}^{1,3}$.}}}
{{\smc Figure~\fci. {\rm The region $H_{7,6}(2)$.}}}
\endinsert

The right hand side of (\ebc) can be evaluated just like in Case 1. Using the same arguments we employed there, we obtain that (see Figure {\fch} for an illustration) 
$$
\spreadlines{3\jot}
\align
\M_{\odot,|}(H_{2a+1,2b}(k_1,\dotsc,k_s))
&=\M_f(D_{a,b,0}^{i_1,\dotsc,i_{a-s}}))\\
&=
\prod_{j=1}^{a-s}{a+b+i_j-1 \choose 2i_j-1}
\prod_{1\leq j<k\leq a-s}\frac{i_k-i_j}{i_j+i_k}.
\tag\ecl
\endalign
$$
It is straightforward to check that the expression given by (\ecj)--(\eck) is equal to the square of the expression on the right hand side of (\ecl). This completes the proof of (\ebc).

%\def\epsfsize#1#2{0.25#1}
%\topinsert
%\centerline{\mypic{3-1x.eps}}
%\centerline{{\smc Figure~\fci. {\rm The region $H_{7,6}(2)$.}}}
%\endinsert

\topinsert
\twoline{\mypic{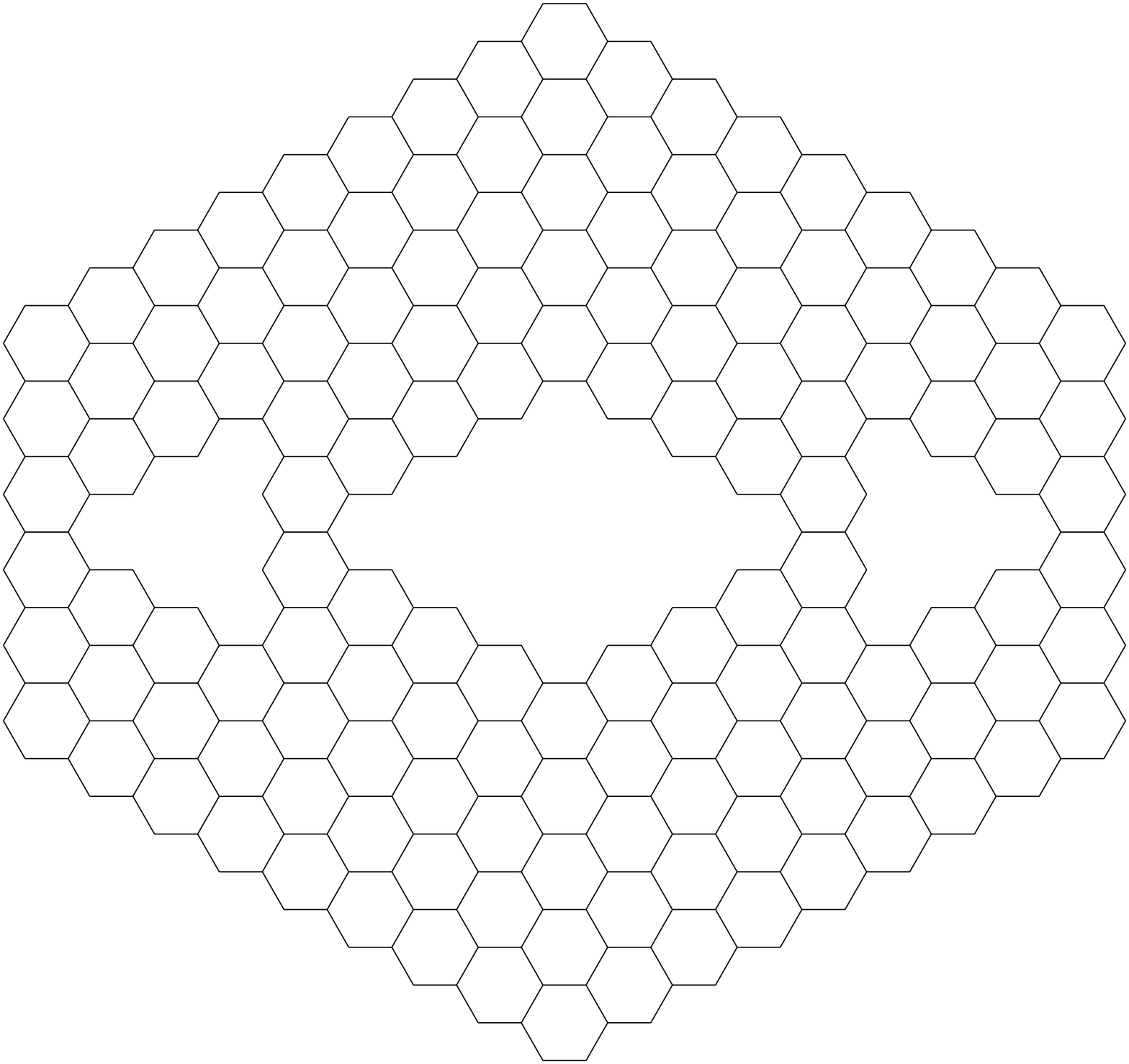}}{\mypic{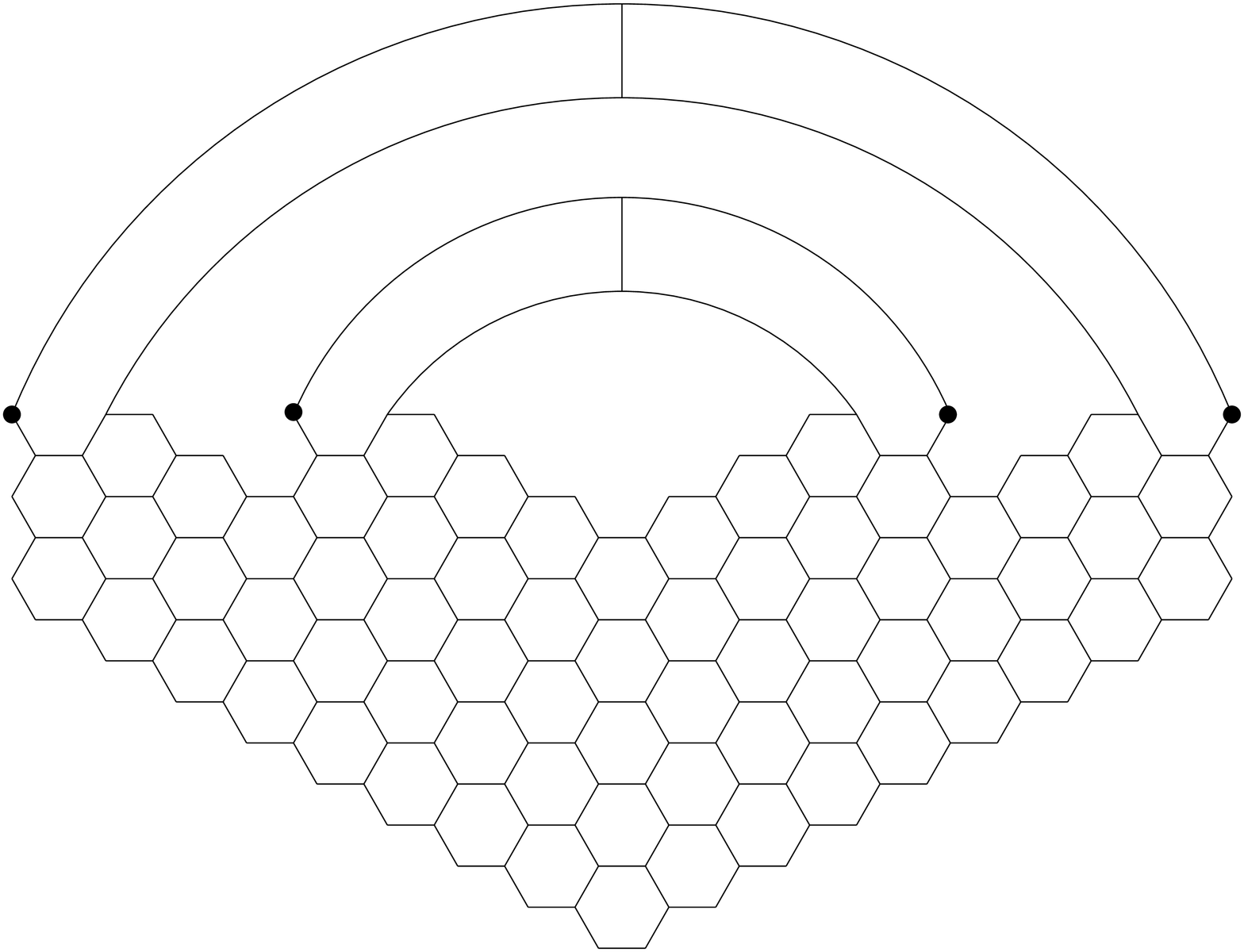}}
\medskip
\centerline{{\smc Figure~{\fcj}. {\rm  The dual graph of $H_{7,6}(2)$, and its quotient by $180^\circ$ rotation.}}}
\endinsert

\topinsert
\twoline{\mypic{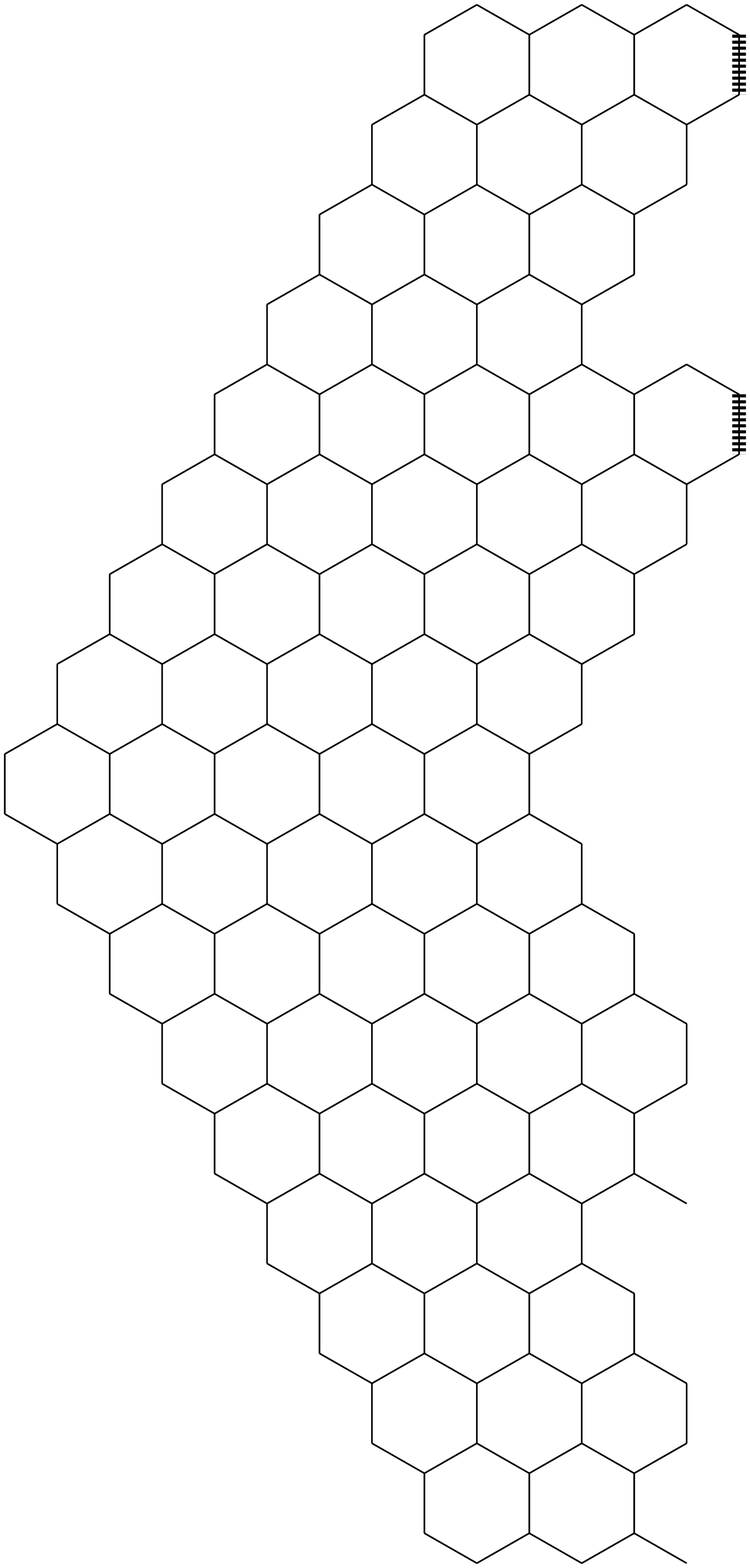}}{\mypic{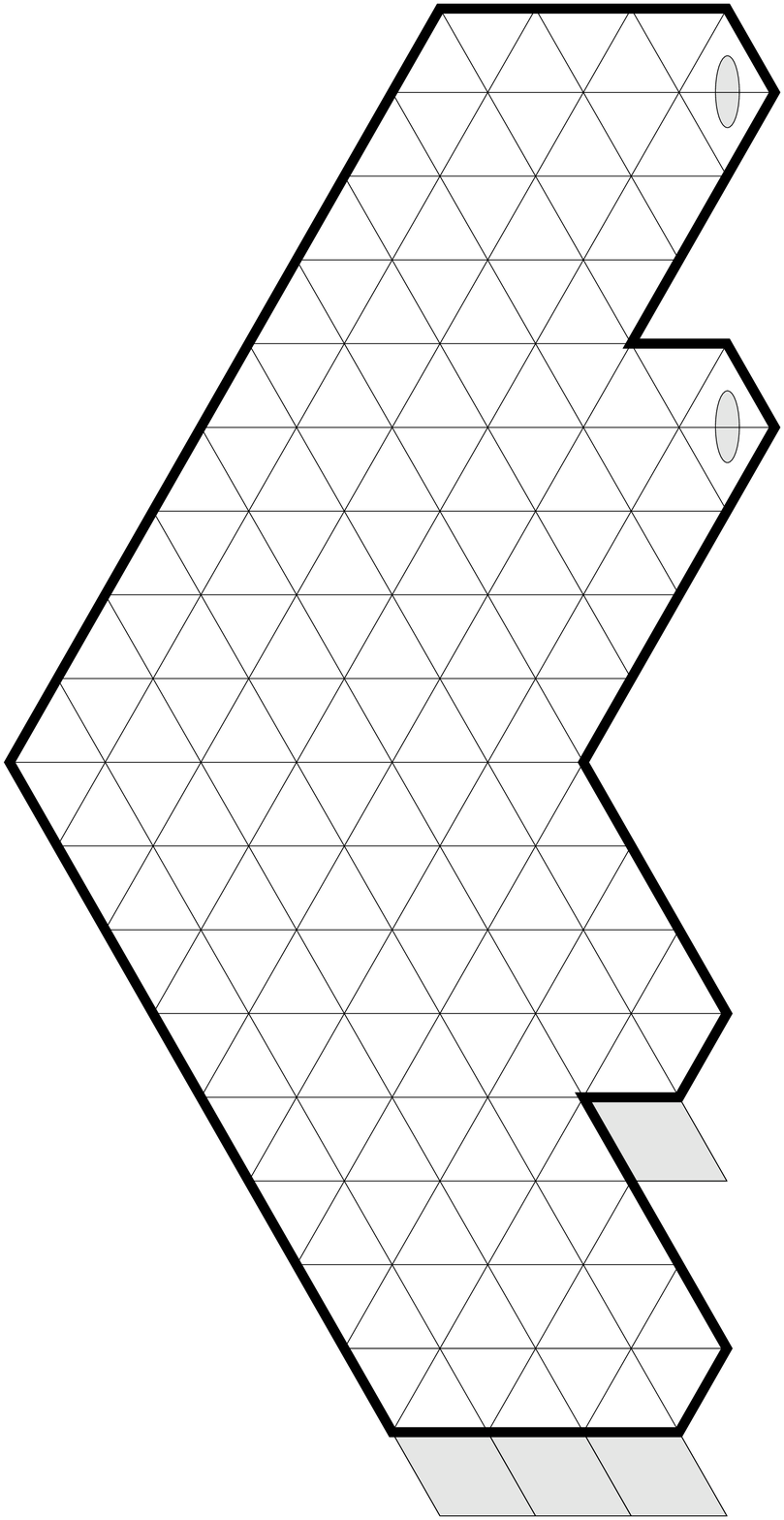}}
\medskip
\centerline{{\smc Figure~{\fck}. {\rm  The graph $K_{7,6}(2)$ resulting by applying the }}}
\centerline{{{\rm factorization theorem to $H_{7,6}(2)$, and its dual region.}}}\endinsert

\topinsert
%\medskip
\centerline{\mypic{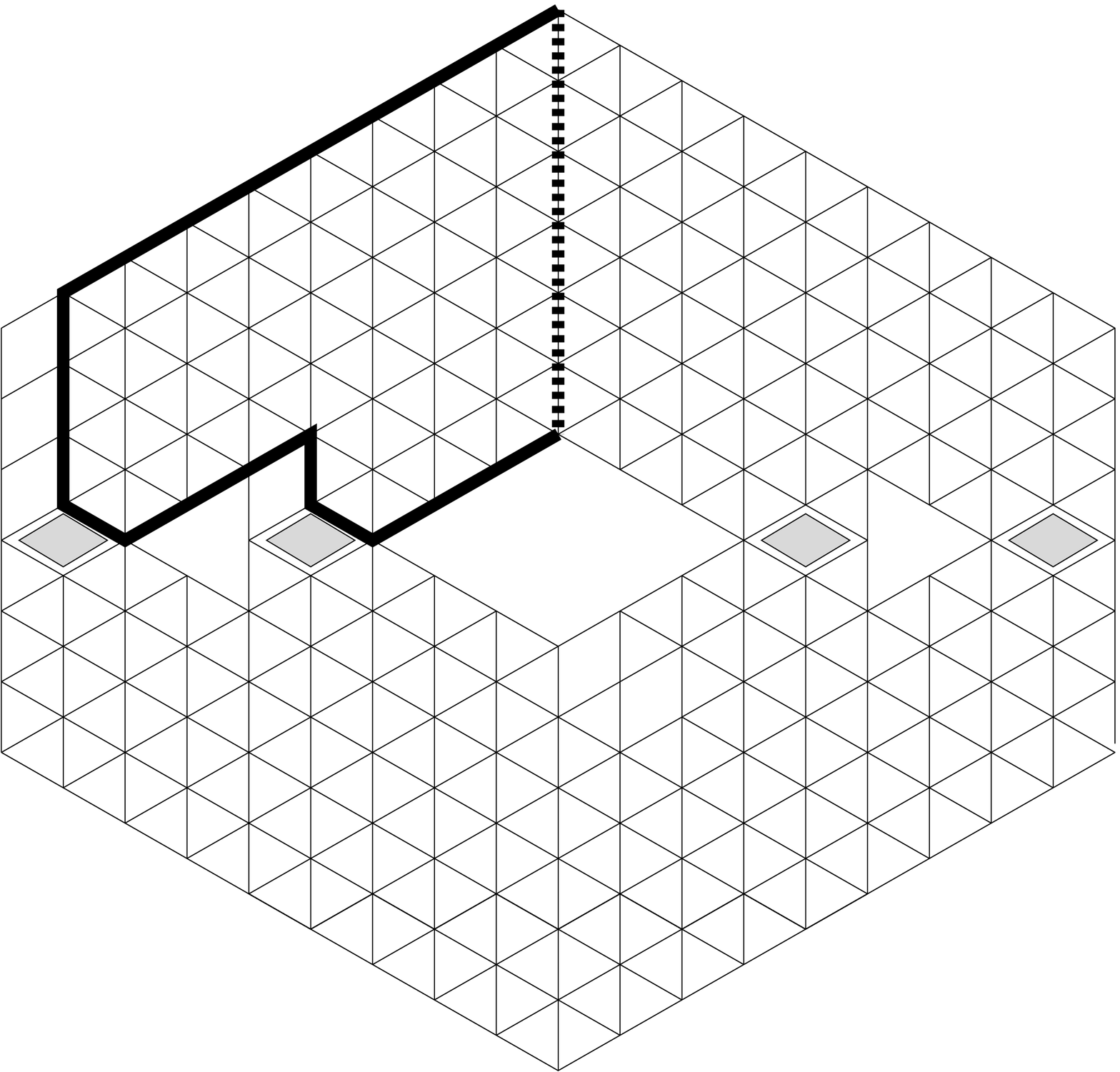}}
\centerline{{\smc Figure~\fcl. {\rm The region $D_{3,3,0}^{1,3}$.}}}
\endinsert
The proof of (\ebd) involves just a very slight modification of the proof of (\ebc) in the case of odd $a$ presented above. The same arguments that led to (\eci) give (see Figures {\fci}--{\fck})
$$
\spreadlines{3\jot}
\align
&
\M(H_{2a+1,2b}(k_1,\dotsc,k_s;2x+1))=
\\
&\ \ \ \ \ \ \ \ \ \ \ \ \ \ \ \ 
2^{a-s}\M(\bar{R}_{[a]\setminus[1,2,\dotsc,x,a-k_1+1,\dotsc,a-k_s+1],[a]\setminus[1,2,\dotsc,x,a-k_1+1,\dotsc,a-k_s+1]}(b)).
\tag\ecm
\endalign
$$
Note that the expression on the right hand side above is a special case of the expression on the right hand side of (\eci). 

Furthermore, the same arguments as in the proof of Case 2 above yield an expression for the right hand side of (\ebd) which is a special case of (\ecl) (see Figure {\fcl}). Therefore, the fact that the two sides of (\ebd) are equal follows from the fact that the expression given by (\ecj)--(\eck) is equal to the square of the expression on the right hand side of (\ecl), i.e. from Case 2 above.

This concludes the proof of Theorem {\tba}.
\epf

\mysec{4. Concluding remarks and open problems}

In this paper we presented a set of four identities involving all ten symmetry classes of plane partitions (see equations (\eaa)--(\ead)), and interpreted them in a uniform way as factorization formulas for the number of perfect matchings of various orbit graphs of honeycomb graphs (see (\eai)--(\eal)). We presented an extension of one of these formulas ((\eaj)), which is the natural counterpart of the extension of another one of them (namely, (\eai)), which we proved in earlier work (see \cite{\fakt}).

Our proof consists of separately evaluating the two sides of the identity, and verifying that they are equal. It would be interesting to have a direct proof of Theorem {\tba}, which shows the equality of the two sides without explicitly evaluating them.

Given the compelling uniformity of equations (\eai)--(\eal), it would be great to find extensions of them that are also uniform. This is accomplished in part (namely, for equations (\eai)--(\eaj)) by \cite{\fakt, Theorem\,2.1} and Theorem {\tba} of this paper. However, the involved generalized regions are not invariant under the extra symmetries required to phrase the analogs of (\eak) and (\eal). It would be interesting to find some other common extension of (\eai)--(\eal) in the same spirit, or at least extensions of (\eak) and (\eal).

\mysec{References}
{\openup 1\jot \frenchspacing\raggedbottom
\roster

%\myref{\anglep}
%  M. Ciucu and I. Fischer, A triangular gap of side 2 in a sea of dimers in a $60^\circ$ angle, {\it J. Phys. A: Math. Theor.} {\bf 45} (2012), 494011.

%\myref{\aanglep}
%  M. Ciucu, A triangular gap of side 2 in a sea of dimers in a $120^\circ$ angle, preprint, 2013.

%\myref{\ff}
% M Ciucu and C. Krattenthaler, A dual of MacMahon's theorem on plane partitions, {\it Proc. Nat. Acad. Sci. USA} {\bf 110} (2013), 4518--4523.

%\myref{\hexnotch}
%  M. Ciucu and I. Fischer, Proof of two conjectures of Ciucu and Krattenthaler on the enumeration of lozenge tilings of hexagons with cut off corners, arxiv preprint arXiv:1309.4640, 2013.

\myref{\AndTSSC}
  G. E. Andrews, Plane Partitions, V: The T.S.S.C.P.P. conjecture, {\it J.
Comb. Theory Ser. A} {\bf 66} (1994), 28--39.

\myref{\FT}
  M. Ciucu, Enumeration of perfect matchings in graphs with reflective symmetry, {\it J. Comb. Theory Ser. A} {\bf 77} (1997), 67--97.

\myref{\csts}
  M. Ciucu, The equivalence between enumerating cyclically symmetric, self-com\-ple\-men\-ta\-ry and totally symmetric, self-complementary plane partitions, {\it J. Comb. Theory, Ser. A} {\bf 86} (1999), 382--389.

\myref{\ppone}
  M. Ciucu, Plane partitions I: A generalization of MacMahon's formula,
{\it Mem. Amer. Math. Soc.} {\bf 178} (2005), no. 839, 107--144.

\myref{\ranglep}
  M. Ciucu, Lozenge tilings with gaps in a $90^\circ$ wedge domain with mixed boundary conditions, {\it Comm. Math. Phys.}, in press.

\myref{\fakt}
  M. Ciucu and C. Krattenthaler, A factorization theorem for rhombus tilings, preprint (2014), arXiv:1403.3323. 

\myref{\MacM}
  P. A. MacMahon, Memoir on the theory of the partition of numbers---Part V. 
Partitions in two-dimensional space, {\it Phil. Trans. R. S.}, 1911, A.

\myref{\DT}
  G. David and C. Tomei, The problem of the calissons, {\it Amer. Math.
Monthly} {\bf 96} (1989), 429--431.

\myref{\KupCSSC}
  G. Kuperberg, Symmetries of plane partitions and the permanent-determinant
method, {\it J. Comb. Theory Ser. A} {\bf 68} (1994), 115--151.

\myref{\StanPP}
  R. P. Stanley, Symmetries of plane partitions, {\it J. Comb. Theory Ser. A}
{\bf 43} (1986), 103--113.

\myref{\StemTS}
  J. R. Stembridge, The enumeration of totally symmetric plane partitions, {\it Adv. in Math.} {\bf 111} (1995), 227--243.

\endroster\par}

\enddocument

\newpage

It turns out that the natural analogs of (\eaa) and (\eac) hold more generally for regions of the type illustrated in Figure {\fab}. More precisely, ... def of regions... Then we have

1. analog of 1.1

2. analog of 1.3.

Indeed, former proved in [fakt]. Main result of this paper is a proof latter.

\newpage

Explanation: How we guessed them:
In this connection, immediate consequence of FT is that $TC|P$. Q: Is ratio something relevant? A: Yes -- it's S!!!

When we noticed this 1 1/2 decades ago, we wanted to see if it's more than a coincidence. Natural to look at hex as 6 60 deg sects. Then we can mod out by Z2, Z3, Z6 --- still get PP numbers on LHS. And RHS has natural analogs. And equality holds!!!!

\newpage

Plane partitions that fit in an $a\times b\times c$ box are well known to be equivalent to lozenge tilings of a hexagon of side-lengths $a$, $b$, $c$, $a$, $b$, $c$ (in cyclic order) on the triangular lattice (for definiteness, we consider that the triangular lattice is drawn in the plane so that one family of parallel lattice lines is vertical). 

Furthermore, the ten symmetry classes of plane partitions correspond to the symmetry classes of lozenge tilings of the corresponding hexagon. More precisely, the generators of the group of symmetries of plane partitions, namely $(i)$ swapping the $x$ and $y$ coordinate axes, $(ii)$ cyclically shifting the axes, and $(iii)$ taking the complement of the diagram of the plane partition in the box enclosing it, correspond to the associated tiling being invariant under reflection across the vertical, rotation by $120$ degrees, and rotation by $180$ degrees, respectively. The formulas of \cite{\StanPP}\cite{\AndTSSC}{\KupCSSC}{\StemTS} give then product expressions for the number of lozenge tilings in each symmetry class.

Translated in this language, equations (\eaa)-(\ead) read
%$$
%\M(H_{a,b})=\M_{-}(H_{a,b})\M_{|}(H_{a,b}),\tag\eaa
%$$
%where $\M(R)$ denotes the number of lozenge tilings of the lattice region $R$, and $\M_{-}(R)$ (resp., $\M_{|}(R)$) denotes the number of those tilings of $R$ that are invariant under reflection across the horizontal (resp. vertical).
%Furthermore, if we require in addition that the tilings involved in (\eaa) be invariant under rotation by 180, 120 or 60 degrees, the equality still holds:
$$
\spreadlines{2\jot}
\align
\M(H_{a,b})&=\M_{-}(H_{a,b})\M_{|}(H_{a,b})\tag\eae\\
\M_{\langle r^3\rangle}(H_{a,b})&=\M_{\langle r^3,-\rangle}(H_{a,b})\M_{\langle r^3,|\rangle}(H_{a,b})\tag\eaf\\
\M_{\langle r^2\rangle}(H_{a,b})&=\M_{\langle r^2,-\rangle}(H_{a,b})\M_{\langle r^2,|\rangle}(H_{a,b})\tag\eag\\
\M_{\langle r\rangle}(H_{a,b})&=\M_{\langle r,-\rangle}(H_{a,b})\M_{\langle r,|\rangle}(H_{a,b}),\tag\eah
\endalign
$$
where $\M(R)$ denotes the number of lozenge tilings of the lattice region $R$, and $\M_{-}(R)$ (resp., $\M_{|}(R)$) denotes the number of those tilings of $R$ that are invariant under reflection across the horizontal (resp. vertical), $r$ is the counterclockwise rotation by $60$ degrees, and $M_{\langle g,h\rangle}(R)$ denotes the number of tilings of the region $R$ that are invariant under group generated by the symmetries $g$ and $h$.

The strikingly uniform nature of identities (\eaa)--(\ead) is inviting of an explanation as to why these identities hold. One way of trying to understand this is to place these identities in a larger context, by finding suitable generalizations of them.

It turns out that the natural analogs of (\eaa) and (\eac) hold more generally for regions of the type illustrated in Figure {\fab}. More precisely, ... def of regions... Then we have

1. analog of 1.1

2. analog of 1.3.

Indeed, former proved in [fakt]. Main result of this paper is an of latter.

\newpage
* noticed $\M(H)=\M_{-}(H)\M_|(H)$

* what about moding out by $\Z^k$, $k=2,3,6$?

* all resulting 3 identities hold!

* interesting feature: all 10 symm classes are involved!

* all 4 are proved --- bec all 10 symm classes are 

* direct proofs? generalizations?

* (i) proved by repr. theory [CK], gen to Schur fn's

* (i) proved algebro-combinatorially, and generalized, in [fakt]

* Main result of this paper: to prove a [fakt]-couterpart gen of (ii) 

* (iii) essentially proved by Stembridge (this was his solution to TS case)

* (iv) proved directly in [the equiv...]

* most natural generalizations of (iii) and (iv) patterned on above common gen of (i) and (ii) don't hold

* open problems: find such gen's for (iii) and (iv); prove directly our gen of (ii) (without separately evaluating LHS, RHS); (completely combinatorial --- bijective? --- proofs of (iii) and (iv); above mentioned proofs are algebro-combinatorial)

* (?) the problem of ``$n$ $60^\circ$ sectors'': Do analogous factorizations hold for any number of sectors different from 1,2,3 and 6? I think we checked $n=4,5$ and it didn't work, maybe also $n=12$.

\enddocument